\pdfoutput=1
\documentclass[5p, times]{elsarticle}

\biboptions{sort&compress}

\usepackage{amsmath, amssymb} 
\usepackage{mathtools}       
\usepackage{bm}              
\usepackage{booktabs}        
\usepackage{graphicx}        
\usepackage{hyperref}        

\usepackage{xurl}            
\usepackage{microtype}       
\usepackage{placeins}

\usepackage{longtable}       
\usepackage{threeparttable}  

\usepackage{pdflscape}       
\usepackage{siunitx}         
\usepackage{chemformula}     
\usepackage{caption}         

\usepackage[official]{eurosym} 
\usepackage{xspace}          
\usepackage{standalone} \usepackage{tikz}
\usepackage{adjustbox}

\usetikzlibrary{shapes.geometric, arrows.meta, positioning, calc, shapes.multipart, fit, backgrounds} 

\usepackage{lineno}
\usepackage{xcolor}

\usepackage[linesnumbered, vlined]{algorithm2e}
\SetKwFor{ParFor}{parfor}{do}{end}%

\SetCommentSty{mycommfont}

\RestyleAlgo{ruled}
\SetKw{GoTo}{go to}

\hypersetup{
    colorlinks=true,
    linkcolor=blue,
    citecolor=blue,
    filecolor=magenta,      
    urlcolor=cyan,
    pdftitle={A Robust Optimization Framework for Energy System Capacity Expansion under Line Outage Uncertainty},
    pdfauthor={Kamran Forghani, et al.},
}

\allowdisplaybreaks 
\newcommand{\setN}{\mathcal{N}}      
\newcommand{\setL}{\mathcal{L}}   
\newcommand{\setLac}{\mathcal{L}^{\text{AC}}}
\newcommand{\setLin}[1]{\mathcal{L}^\text{in}_{#1}}      
\newcommand{\setLout}[1]{\mathcal{L}^\text{out}_{#1}}      
\newcommand{\setT}{\mathcal{T}}      
\newcommand{\setTzero}{\mathcal{T}^0} 
\newcommand{\setGgen}{\mathcal{G}}
\newcommand{\setO}{\Omega}  

\newcommand{\setX}{\mathcal{X}}

\newcommand{\load}[2]{D_{#1, #2}}              
\newcommand{\cf}[3]{\mathrm{cf}_{#1, #2, #3}}     
\newcommand{\maxcap}[2]{\overline{C}_{#1, #2}} 
\newcommand{\mincap}[2]{\underline{C}_{#1, #2}}
\newcommand{\maxcaptrans}[1]{\overline{F}_{#1}} 
\newcommand{\mincaptrans}[1]{\underline{F}_{#1}} 
\newcommand{\ANNUALCOST}[1]{c^{\text{inv}}_{#1}} 
\newcommand{\ANNUALCOSTLine}[1]{c^{\text{tra}}_{#1}} 
\newcommand{\VAROMCOST}[1]{c^{\text{var}}_{#1}}   
\newcommand{\EMISSIONFACTOR}[1]{e_{#1}}             
\newcommand{\CAPACITYCREDIT}[1]{\nu_{#1}}           
\newcommand{\COSTSHED}{c^{\text{shed}}} 
\newcommand{\COSTLOSS}{c^{\text{loss}}} 
\newcommand{\MAXDEMANDCURTAILMENTRATIO}{\delta^{\text{shed}}} 
\newcommand{\EFFICIENCYCHARGE}{\eta^{+}}
\newcommand{\EFFICIENCYDISCHARGE}{\eta^{-}}
\newcommand{\Bline}[1]{B_{#1}} 
\newcommand{\Lline}[1]{L_{#1}}
\newcommand{\CARBONCAPINTENSITY}{E^{\text{CO}_2}}    
\newcommand{\HOURSPERYEAR}{H_{\text{yr}}}          
\newcommand{\linestatus}[2]{\xi_{#1, #2}} 
\newcommand{\LINESTATUS}{\bm{\xi}} 
\newcommand{\RESERVEMARGIN}{\Gamma}  

\newcommand{\TIMESTEPDURATION}{\Delta t}
\newcommand{\OPWEIGHT}{w}   

\newcommand{\COSTX}{\bm{c}}
\newcommand{\COSTy}{\bm{q}}
\newcommand{\bx}{\bm{b}}
\newcommand{\Ax}{\bm{A}}
\newcommand{\Wy}{\bm{W}}
\newcommand{\hy}{\bm{h}}
\newcommand{\Ty}[1]{\bm{T}(#1)}

\newcommand{\Captech}[2]{C_{#1, #2}}                
\newcommand{\Captrans}[1]{F_{#1}}
\newcommand{\CapEl}[1]{C^{\text{el}}_{#1}}          
\newcommand{\CapFC}[1]{C^{\text{fc}}_{#1}}          
\newcommand{\CapHS}[1]{C^{\text{H}_2}_{#1}}         
\newcommand{\INVESTMENT}{\bm{x}}
\newcommand{\OPERATION}{\bm{y}}
\newcommand{\Qx}[2]{Q(#1, #2)}
\newcommand{\TOC}{\tau_{\max}} 

\newcommand{\Pgen}[3]{p_{#1, #2, #3}}                
\newcommand{\Estor}[2]{E_{#1, #2}} 
\newcommand{\Cin}[2]{I^+_{#1, #2}} 
\newcommand{\Cout}[2]{I^-_{#1, #2}} 
\newcommand{\LineFlowFwd}[2]{f^{+}_{#1, #2}} 
\newcommand{\LineFlowBwd}[2]{f^{-}_{#1, #2}} 
\newcommand{\ThetaVar}[2]{\theta_{#1, #2}} 
\newcommand{\Pshed}[2]{p^{\text{shed}}_{#1, #2}} 
\newcommand{\Ploss}[2]{p^{\text{loss}}_{#1, #2}} 

\newcommand{\dualGen}[3]{\pi^{\text{gen}}_{#1,#2,#3}}       
\newcommand{\dualEl}[2]{\pi^{\text{el}}_{#1,#2}}           
\newcommand{\dualFC}[2]{\pi^{\text{fc}}_{#1,#2}}           
\newcommand{\dualHS}[2]{\pi^{\text{H}_2}_{#1,#2}}          
\newcommand{\dualUns}[2]{\pi^{\text{uns}}_{#1,#2}}         
\newcommand{\dualLine}[2]{\pi^{\text{line}}_{#1,#2}}        
\newcommand{\dualNode}[2]{\pi^{\text{node}}_{#1,#2}}        
\newcommand{\dualShed}[2]{\pi^{\text{shed}}_{#1,#2}}        
\newcommand{\dualCO}{\pi^{\text{CO2}}}                     
\DeclareSIUnit{\TonCOtwo}{\text{t}_{\ch{CO2}}}
\DeclareSIUnit{\GramCOtwo}{\text{g}_{\ch{CO2}}}
\DeclareSIUnit{\MTonCOtwo}{\text{Mt}_{\ch{CO2}}}
\DeclareSIUnit{\KgCOtwo}{\text{kg}_{\ch{CO2}}}
\DeclareSIUnit{\Currency}{\text{\euro}\xspace}
\DeclareSIUnit{\year}{\text{y}}
\DeclareSIUnit{\MW}{\text{MW}}
\DeclareSIUnit{\GW}{\text{GW}}
\DeclareSIUnit{\KW}{\text{kW}}
\DeclareSIUnit{\kWh}{\text{kWh}}
\DeclareSIUnit{\MWh}{\text{MWh}}
\DeclareSIUnit{\GWh}{\text{GWh}}
\DeclareSIUnit{\TWh}{\text{TWh}}
\DeclareSIUnit{\MCurrency}{\text{M\euro}\xspace}

\newcommand{\VoRCost}{VoR\textsuperscript{Cost}\xspace}
\newcommand{\VoRLoL}{VoR\textsuperscript{LoL}\xspace}

\newcommand{\noIterations}{I^{\max}}
\newcommand{\noRandomScenarios}{N}
\newcommand{\noMemoryScenarios}{M}
\newcommand{\noSimScenarios}{N^{\text{sim}}} 
\newcommand{\ConfidenceLevel}{\alpha} 
\newcommand{\SampledScenarios}{\mathcal{S}} 
\newcommand{\VaR}{\mathrm{VaR}}

\newcommand{\LB}{LB}
\newcommand{\UB}{UB}
\newcommand{\Gap}{\text{Gap}}
\newcommand{\RandomScenarios}{\mathcal{A}}
\newcommand{\MemoryScenarios}{\mathcal{M}}
\newcommand{\DualVariables}{\bm{\pi}}
\newcommand{\TC}{TC}

\newcommand{\MTTR}{\text{MTTR}}
\newcommand{\FailRateBase}{\lambda_{0}}
\newcommand{\ProbRep}[1]{\rho^{\text{rep}}_{#1}}
\newcommand{\ProbFail}[1]{\rho_{#1}}
\newcommand{\FailRate}[1]{\lambda_{#1}}
\newcommand{\RepairRate}{\mu}

\newcommand{\ProbInit}[1]{\rho^{0}_{#1}}
\newcommand{\Indicator}[1]{\mathbb{I}(#1)}
\newcommand{\RandU}[2]{\mathrm{Uniform}(#1, #2)}

\begin{document}

\begin{frontmatter}

\title{A Robust Optimization Framework for Energy System Capacity
Expansion under Line Outage Uncertainty}

\author[chalmers]{Kamran Forghani\corref{cor1}}
\cortext[cor1]{Corresponding author}
\ead{kamranf@chalmers.se}

\author[chalmers]{Lina Reichenberg}
\author[chalmers]{Fredrik Hedenus}
\author[Berlin]{Iegor Riepin}

\affiliation[chalmers]{
            organization={Department of Environmental and Energy Sciences, Chalmers University of Technology},
            city={Göteborg},
            country={Sweden}}

\affiliation[Berlin]{
            organization={Department of Digital Transformation in Energy Systems, TU Berlin},
            country={Germany}}

\begin{abstract}
    Capacity expansion models are widely used as a tool in climate and energy policy analysis to find optimal investment decisions in generation, storage, and transmission. As the share of renewable generation increases in energy systems, these models rely more on transmission lines to balance supply and demand across regions. Transmission lines, however, are at risk of failure due to aging, extreme weather, or accidents, and such failures undermine this load balancing. Models that represent outages as time-dependent events, which capture when and for how long lines are out of service, come mainly from electrical engineering and are applied to small-scale systems, for example, in reliability assessment and short-term operational planning. These models, however, do not scale easily to inter-regional capacity expansion, as the optimization problem quickly becomes computationally intractable. In this paper, we develop a capacity expansion model in which the risk of line failures shapes the investments in all parts of the system. We formulate it as a two-stage, scenario-based robust optimization problem that samples time-varying outage scenarios from a Markov model and hedges the investments against the most severe ones. To solve the model, we propose a memory-augmented Benders decomposition algorithm combined with Monte Carlo simulation that can handle thousands of scenarios. The algorithm then evaluates the selected plan with an out-of-sample simulation to bound the Value-at-Risk (VaR) of its system cost. In a 25-node case study of Germany and its neighboring countries, the gap between the lower bound and the simulated VaR bound is below 1\%. We compare the robust plans with a deterministic plan, which assumes all lines are always available. We also compare them with a reserve-margin heuristic, which requires local firm capacity to be above peak demand. Like the heuristic, the robust plans nearly eliminate loss of load in 1000 simulated outage scenarios. The mean system cost of the robust plans is within 0.2\% of the deterministic plan, but about 4.5\% below the heuristic. Instead of adding more generation capacity, the robust approach mainly changes the type and location of capacity, for example, by spreading storage across the network.
\end{abstract}


\end{frontmatter}

\section{Introduction} \label{sec:Introduction}

Energy system optimization models are important tools for exploring and evaluating pathways toward a secure, low-carbon energy system \citep{3_yue_2018_review}. Among these models, the Generation and Transmission Expansion Problem (G\&TEP), also known as the capacity expansion problem, has been widely used to identify least-cost transition strategies \citep{7_in_3_giannakidis2015informing}. The G\&TEP typically determines long-term investments in generation (e.g., wind, solar, and thermal plants), storage, and transmission that minimize the system cost while meeting operational constraints. The results of such models can help evaluate decarbonization pathways and their economic and environmental effects on future energy systems \citep{7_in_3_giannakidis2015informing, 3_yue_2018_review}.

Most capacity expansion models for inter-regional planning apply a deterministic approach by assuming all parameters are known in advance \citep{3_yue_2018_review, 30_roald2023ppower}. Many parameters of real energy systems, however, such as load, resource availability, and technology costs, are uncertain. Several approaches have been developed in the literature to deal with this uncertainty, and stochastic programming and robust optimization are the most common ones \citep{3_yue_2018_review, 30_roald2023ppower}. In stochastic programming, uncertainty is represented explicitly through scenarios. It is used in established models such as EMPIRE \citep{40_Backe_2022_EMPIRE} and stochastic extensions of SWITCH \citep{41_Johnston_2019_Switch}, and it has recently been added to PyPSA \citep{pypsa_stochastic_optimization_2023}. Further studies on the applications of stochastic and robust methods for representing uncertainty in resource availability, demand, and technology costs in energy system models can be found in \citep{44_moret2020decision, 45_verastegui2019adaptive, 2_mantegna_2026_uncertainty, 43_Kang_2023_stochastic, forghani2025modeling}. To our knowledge, however, none of these inter-regional G\&TEP formulations considers uncertainty in the \emph{availability} of transmission lines, which is the main focus of this paper.

As geographically dispersed renewable generation becomes more prevalent, transmission networks become more important for balancing supply and demand across regions \citep{8_Wang2022_dist, 21_Wen2022_robust}. Transmission lines, however, can fail for various reasons such as aging infrastructure, extreme weather, or accidents \citep{5_Moreira2015_ARO, 31_tapia_2021_robust, 30_roald2023ppower}. Energy system planning should therefore include such contingencies, so that investment decisions remain robust when failures occur.

\subsection{Related Work and Research Gap} \label{sec:LiteratureReview}

To keep an energy system operational when components fail, the $N$-$k$ security criterion was proposed. It requires that the system remain operational after the loss of any $k$ or fewer components \citep{5_Moreira2015_ARO}. Its most common special case, the $N$-1 criterion, is a standard reliability requirement of the North American Electric Reliability Corporation \citep{nerc_standards} and is widely applied in practice and in planning models because of its computational simplicity. In inter-regional energy system optimization models, $N$-1 security is usually enforced through deterministic security constraints. Linear optimal power flow (LOPF) formulations, for example, can include the $N$-1 criterion, as implemented in PyPSA \citep{42_Gazafroudi_2022_topology, pypsa_scigrid_sclopf}. Studies of smaller electrical systems also enforce the $N$-1 criterion deterministically, through predefined contingency sets or hard constraints \citep{14_oboudi2025reliability, 20_Paula2020, 13_Peng2023, 16_MajidiQadikolai2016, 12_Mehrtash2020}.

Growing exposure to extreme weather and the need for more secure energy systems lead researchers to consider several simultaneous failures, as in the general $N$-$k$ criterion. This generalization, however, scales poorly beyond $N$-1. The reason is that the number of contingencies grows quickly, so they must either be reduced considerably or fully enumerated.\footnote{The number of scenarios grows quickly with $k$. In a system with 20 lines, for example, the $N$-1 criterion requires $\binom{20}{1} = 20$ scenarios, while $N$-2 and $N$-3 add $\binom{20}{2} = 190$ and $\binom{20}{3} = 1140$ scenarios. In large energy systems with hundreds of components, evaluating all $N$-$k$ contingencies therefore becomes computationally intractable.} Monte Carlo sampling combined with scenario reduction has been used to keep the number of contingency scenarios reasonable \citep{39_su_2024_stochastic}. Robust optimization, combined with decomposition methods, has also been applied to contingency-constrained expansion planning. These approaches generally search for solutions that remain feasible under the worst realization of uncertain parameters within an uncertainty set. In most of these studies, an uncertainty budget is imposed to limit the number of simultaneous failures \citep{5_Moreira2015_ARO, 9_Wu2018_robust, 6_Liu2022_distributed}. Robust optimization techniques are often criticized for their conservatism because failures can occur simultaneously under adverse operating conditions; combinations that might never occur in reality but can lead to overinvestment. Distributionally robust optimization can reduce this conservatism by modeling ambiguity in the probability distribution of contingencies \citep{8_Wang2022_dist, 4_Yuan2022_TEP} or in renewable generation \citep{18_Zhang2021_data}, which gives a better trade-off between robustness and cost.

Table~\ref{tab:lit_review} summarizes the studies in our review that address transmission contingencies in a way closely related to G\&TEP or the transmission expansion problem (TEP). We use these studies as methodological inspiration rather than as directly comparable work. Most of them address different research questions and are tested on small, standard test networks, not on inter-regional energy systems.

\begin{table*}[t]
    \centering
    \caption{Overview of transmission and capacity expansion models with transmission line contingencies.} \label{tab:lit_review}
    \resizebox{\textwidth}{!}{
    \begin{threeparttable}
        \begin{tabular}{@{} l l l l l l l l l @{}}
            \toprule
            \textbf{Paper} & \textbf{Problem} & \textbf{Criterion} & \textbf{Contingency Representation} & \textbf{Other Uncertainty} & \textbf{Modeling Approach} & \textbf{Outage Duration} & \textbf{Storage} & \textbf{Solution Method} \\
            \midrule
            \citep{5_Moreira2015_ARO}        & TEP         & $N$-$k$\tnote{a}  & Uncertainty budget                 & --                  & RO (worst-case)                 & No  & No  & BD + C\&CG        \\
            \citep{16_MajidiQadikolai2016}   & TEP         & $N$-1             & Screened list (VCL/CII)            & --                  & Deterministic                   & No  & No  & MILP (iterative)  \\
            \citep{9_Wu2018_robust}          & TEP         & $N$-$k$           & Uncertainty budget                 & Load, RES           & RO (worst-case)                 & No  & No  & BD                \\
            \citep{20_Paula2020}             & TEP         & $N$-1             & Enumeration                        & Wind                & RO (scenario-based)             & No  & No  & BD                \\
            \citep{12_Mehrtash2020}          & G\&TEP      & $N$-1             & Enumeration with propagation       & --                  & Deterministic + risk            & No  & No  & MISOCP            \\
            \citep{18_Zhang2021_data}        & TEP         & $N$-$k$\tnote{a}  & Uncertainty budget                 & RES                 & DRO (worst-case expectation)    & No  & No  & C\&CG             \\
            \citep{6_Liu2022_distributed}    & DTEP        & $N$-$k$           & Uncertainty budget                 & Load, RES           & RO (worst-case) + stochastic    & No  & Yes & C\&CG             \\
            \citep{8_Wang2022_dist}          & TEP         & $N$-$k$\tnote{a}  & Scenarios + uncertainty budget     & --                  & DRO (worst-case expectation)    & No  & No  & BD                \\
            \citep{21_Wen2022_robust}        & TEP         & $N$-1             & Enumeration (LODF)                 & Load, RES           & RO (worst-case)                 & No  & No  & C\&CG             \\
            \citep{4_Yuan2022_TEP}           & TEP         & $N$-$k$\tnote{a}  & Prob.\ distribution + budget       & Load, wind          & DRO (worst-case expectation)    & No  & No  & BD                \\
            \citep{13_Peng2023}              & G\&TEP      & $N$-1             & Enumeration                        & Load, RES           & RO (worst-case)                 & No  & No  & Metaheuristic     \\
            \citep{39_su_2024_stochastic}    & T\&PEP      & $N$-$k$\tnote{b}  & Monte Carlo scenarios + reduction  & --                  & Stochastic                      & No  & No  & MINLP             \\
            \citep{14_oboudi2025reliability} & TEP         & $N$-1             & Reliability constraints            & Load, RES           & Stochastic                      & No  & No  & Metaheuristic     \\
            \midrule
            \textbf{This paper}              & G\&TEP      & $N$-$k$\tnote{c}  & Monte Carlo scenarios (Markov)     & --                  & RO (sampled worst-case)         & Yes & Yes & BD                \\
            \bottomrule
        \end{tabular}
        \begin{tablenotes}[flushleft]
            \item[a] Outages of both generators and transmission lines.
            \item[b] Geographically correlated failures caused by extreme weather.
            \item[c] Any number of lines can be out of service at the same time, and outages last over several time steps.
            \item \textbf{Abbreviations:}
            \textbf{BD}: Benders Decomposition;
            \textbf{C\&CG}: Column-and-Constraint Generation;
            \textbf{DRO}: Distributionally Robust Optimization;
            \textbf{DTEP}: Distribution and Transmission Expansion Problem;
            \textbf{G\&TEP}: Generation and Transmission Expansion Problem;
            \textbf{LODF}: Line Outage Distribution Factor;
            \textbf{MILP}: Mixed-Integer Linear Programming;
            \textbf{MINLP}: Mixed-Integer Nonlinear Programming;
            \textbf{MISOCP}: Mixed-Integer Second-Order Cone Programming;
            \textbf{$N$-1 / $N$-$k$}: Security criterion representing the loss of 1 or $k$ system components;
            \textbf{RES}: Renewable Energy Sources;
            \textbf{RO}: Robust Optimization;
            \textbf{T\&PEP}: Transmission and Pipeline Expansion Problem;
            \textbf{TEP}: Transmission Expansion Problem;
            \textbf{VCL/CII}: Variable Contingency List / Contingency Identification Index.
        \end{tablenotes}
    \end{threeparttable}}
\end{table*}

Our review of the existing models leads to two observations that motivate this work. The first and main one is that none of the reviewed G\&TEP or TEP studies represent transmission line outages as time-dependent events. As seen in Table~\ref{tab:lit_review}, all of these studies treat failures as static snapshots, independent of when they occur. This simplification is justified for traditional fossil-based energy systems, in which dispatchable generation can respond flexibly to outages and planning models are typically operated at low temporal resolution. However, this simplification is much less suitable for systems with high shares of wind and solar, in which electricity output varies considerably across hours and seasons. In such systems, the consequences of an outage, and thus the value of investments that mitigate its risk, depend on \emph{when} it occurs and \emph{how long} it lasts, not only on \emph{whether} it occurs. A few studies do, in fact, model the duration of outages over time, for example in microgrid scheduling \citep{24_yuanxiong_islanding, 27_Haifeng_bilevel} or in reliability assessment with sequential Monte Carlo simulation \citep{7_Beyza_characterising}. But these approaches work at the scale of a single node or a single system. Therefore, they are not readily applicable to long-term, inter-regional G\&TEP with a chronological representation of demand and renewable generation. Using them in inter-regional planning would require combining their detailed, single-system treatment of outage duration with the multi-node, chronological structure of capacity expansion models.

The second observation is that most of the robust optimization studies reviewed in this study rely on an \emph{uncertainty budget}, which sets how conservative the resulting plan should be. A larger budget results in more robustness, but at a higher cost. The uncertainty budget is typically chosen arbitrarily, often without any statistical interpretation. Most of these models also provide no direct feedback on how an investment plan performs under uncertainty. Instead, performance must be assessed after optimization, usually by Monte Carlo simulation (e.g., \citep{7_Beyza_characterising}), which evaluates a given investment plan but does not influence the optimization. As a result, the models may return overly conservative or costly solutions without exploring alternatives that offer a better trade-off between cost and robustness. Integrating simulation directly into the optimization could address this limitation. By using statistically plausible outage realizations, rather than an externally chosen budget, the optimization could determine the scenarios against which a plan is hedged. This approach would also enable quantification of the risk associated with each candidate plan during optimization, rather than only after completion.

\subsection{Research Questions and Contributions}

The underlying question of this paper is how investments in generation, storage, and transmission should be planned when it is known that transmission lines sometimes fail. To answer it and to address the gaps identified in the literature, we pose the following research questions:
\begin{enumerate}
    \item \label{rq:integration} How can the timing and duration of transmission line outages be included in a robust, inter-regional generation and transmission expansion model for long-term planning?
    \item \label{rq:solvability} How can this robust model be solved efficiently for thousands of outage scenarios?
    \item \label{rq:comparison} How do the robust plans compare with the deterministic and reserve-margin plans in terms of cost and robustness?
\end{enumerate}

The paper makes the following contributions. For RQ\ref{rq:integration}, we formulate a two-stage, scenario-based robust optimization problem (Sections~\ref{sec:ProblemStatement}--\ref{sec:meth_simulation}). The first stage determines long-term investments in generation, storage, and transmission. The second stage optimizes system operation under outage scenarios sampled from a chronological outage model that captures both outage timing and duration. The investment plan is hedged against the most severe sampled scenarios, which gives a tractable approximation of a Value-at-Risk (VaR) criterion on system cost. For RQ\ref{rq:solvability}, we propose a Benders decomposition algorithm combined with Monte Carlo simulation to efficiently explore thousands of scenarios (Section~\ref{sec:benders_decomposition}). In each iteration, the algorithm samples new scenarios and keeps a memory of the most critical ones. It uses Cantelli's inequality to bound the VaR of the system cost of each candidate plan. The sub-problems, one for each outage realization, are solved in parallel and provide optimality cuts for the master problem. The selected plan is finally evaluated on a large set of new, independently simulated scenarios. For RQ\ref{rq:comparison}, we apply the framework to a 25-node, inter-regional case study of Germany and its neighboring countries (Section~\ref{sec:case_study}). We then compare the robust plans with deterministic planning and a reserve-margin heuristic in terms of cost and operational risk, measured by loss of load (Section~\ref{sec:results}).
\section{Methods} \label{sec:methods}

This section presents the methods of the paper. Section~\ref{sec:ProblemStatement} states the planning problem and the modeling assumptions, and Section~\ref{sec:robust_model} formulates the problem as a scenario-based robust optimization model. Section~\ref{sec:meth_simulation} describes the Markov-chain simulation that generates the outage scenarios, and Section~\ref{sec:benders_decomposition} presents the memory-augmented Benders decomposition algorithm that solves the model.

\subsection{Problem Statement and Assumptions} \label{sec:ProblemStatement}

We consider a G\&TEP with a set of network nodes $\setN$, transmission lines $\setL$, a portfolio of generation technologies $\setGgen$, and hydrogen storage. The goal is to determine an optimal capacity investment plan, specifying the generation capacity $\Captech{n}{g}$ and hydrogen storage capacity at each node, as well as the capacity $\Captrans{l}$ of each transmission line. Hydrogen storage consists of three components: electrolyzers that convert electricity into hydrogen (charging power $\CapEl{n}$), hydrogen storage (energy capacity $\CapHS{n}$), and fuel cells that convert hydrogen back into electricity (discharging power $\CapFC{n}$). We propose a robust optimization approach to account for transmission line outages, characterized not only by their occurrence probabilities but also by their durations. The objective is to minimize the system cost while maintaining operational feasibility under adverse transmission availability conditions. The complete list of the mathematical notation used in this paper is provided in the nomenclature in~\ref{app:nomenclature}.

To keep the model tractable while retaining the main physical characteristics of electricity networks and isolating the effect of transmission outages, we make the following assumptions. The implications of these assumptions for the results are discussed in Section~\ref{sec:limitations}.

\begin{enumerate}
    \item \label{ass:single_year} \label{ass:continuous} \textbf{Single-Year, Continuous Capacity Expansion.}
    Capacities are optimized for a single target year with annualized investment costs. All capacities are continuous decision variables, as is common in inter-regional capacity expansion models, to avoid the computational burden of integer variables \citep{horsch2018pypsa, 38_Neumann2022_assessments, Grochowicz2023, 40_Backe_2022_EMPIRE, 41_Johnston_2019_Switch}.

    \item \label{ass:uncertainty} \label{ass:foresight} \textbf{Transmission Outages as the Only Uncertainty.}
    To isolate the effect of transmission outages, fixed time series of load $\load{n}{t}$ and renewable capacity factors $\cf{n}{g}{t}$ are used, with the assumption that generation and storage units are always available. Although demand and weather uncertainty are important for long-term planning \citep{22_Zhang2018_robust, 2_mantegna_2026_uncertainty}, they are beyond the scope of this study. For each outage scenario, operation is optimized with perfect foresight of that scenario.

    \item \label{ass:corridor} \label{ass:outage_process} \textbf{Corridor-Level, Independent Markovian Outages.}
    Each node in the system is an aggregated region, and each element of $\setL$ is an inter-nodal corridor that may aggregate several parallel circuits. Each corridor is either fully available or fully out of service ($\linestatus{l}{t} \in \{0,1\}$). The status of each corridor follows a two-state Markov process with a constant failure rate proportional to its length \citep{19_yang2019failure} and a constant repair rate  (Section~\ref{sec:meth_simulation}). Failures are independent across corridors and do not vary with weather. Note that hereafter, we refer to each corridor $l \in \setL$ as line~$l$.

    \item \label{ass:dcpf} \textbf{Linearized Power Flow.}
    Power flows on AC lines follow the linearized DC power-flow approximation \citep{38_Neumann2022_assessments}. Flows on DC links are fully controllable and are limited only by their capacity and availability (transport model). The approximation neglects reactive power, voltage magnitudes, and transmission losses, but it captures congestion and loop flows and keeps the model linear.

    \item \label{ass:shedding} \textbf{Load Shedding and Loss of Load.}
    Voluntary load shedding $\Pshed{n}{t}$ is allowed at a cost $\COSTSHED$ and is limited to a share $\MAXDEMANDCURTAILMENTRATIO$ of demand. Involuntary loss of load $\Ploss{n}{t}$ is also allowed to keep every outage scenario feasible, but it is heavily penalized.
\end{enumerate}

\subsection{Scenario-Based Robust Optimization Model} \label{sec:robust_model}

Classical worst-case robust optimization enforces feasibility under every combination of failures in a given uncertainty set, which often leads to overly conservative and costly designs \citep{30_roald2023ppower, 8_Wang2022_dist, 18_Zhang2021_data}. We instead hedge only against statistically plausible outage realizations sampled from the outage model. The goal is to limit the risk of high operational costs under such outages. Let $\INVESTMENT$ denote the vector of first-stage investment decisions ($\Captech{n}{g}$, $\CapEl{n}$, $\CapFC{n}$, $\CapHS{n}$, and $\Captrans{l}$). These decisions are constrained by the deterministic polyhedral set $\setX = \left\{ \INVESTMENT \ge \bm{0} \mid \Ax \INVESTMENT \le \bx \right\}$, which corresponds to constraint sets~\eqref{eq:inv_res_potential}--\eqref{eq:inv_h2_fc_el}, and are made before the line outages are known. Uncertainty is represented by $\LINESTATUS \in \setO$, where $\linestatus{l}{t} = 1$ if line $l$ is available in time step $t$, and $0$ otherwise. The set $\setO = \{0,1\}^{|\setL| \times |\setT|}$ contains all possible line status realizations over the time horizon.

For a fixed investment decision $\bar{\INVESTMENT}$ and realization $\LINESTATUS$, the second-stage operational cost $\Qx{\bar{\INVESTMENT}}{\LINESTATUS}$ is the minimum operating cost under the transmission outages implied by $\LINESTATUS$. This cost corresponds to the optimal value of the following linear program.
\begin{align}
    \Qx{\bar{\INVESTMENT}}{\LINESTATUS} = \min_{\OPERATION \ge \bm{0}} \quad & \COSTy^T \OPERATION \nonumber \\
    \text{s.t.} \quad & \Wy(\LINESTATUS) \OPERATION \le \hy - \Ty{\LINESTATUS} \bar{\INVESTMENT}, \label{problem:compact_subproblem}
\end{align}
where $\OPERATION$ denotes the second-stage operational variables, including dispatch ($\Pgen{n}{g}{t}$), hydrogen storage operation ($\Cin{n}{t}$, $\Cout{n}{t}$, and $\Estor{n}{t}$), flows ($\LineFlowFwd{l}{t}$ and $\LineFlowBwd{l}{t}$), and load shedding ($\Pshed{n}{t}$). The vector $\hy$ represents the right-hand sides of the operational constraints, while $\Ty{\LINESTATUS}$ maps the line status realization to the corresponding transmission capacity limits. The matrix $\Wy(\LINESTATUS)$ contains the coefficients of the operational variables and depends on the line status, since $\linestatus{l}{t}$ activates or deactivates the power-flow equation for each AC line (see Constraint~\eqref{eq:op_dc_flow}).

To model uncertainty, consider $(\setO, 2^{\setO}, \mathbb{P})$ as a probability space, where $\mathbb{P}: 2^{\setO} \rightarrow [0,1]$ denotes the probability measure induced by the outage process described in Assumption~\ref{ass:outage_process}. The proposed robust framework is motivated by the Value-at-Risk (VaR) of the operational cost at a probability level $\alpha \in (0,1)$, defined as $\VaR_{\alpha}\big(\Qx{\INVESTMENT}{\LINESTATUS}\big) = \min\{\tau \in \mathbb{R} : \mathbb{P}(\Qx{\INVESTMENT}{\LINESTATUS} \le \tau) \ge \alpha\}$ \citep{rockafellar2002cvar}. This value represents the lowest cost threshold that is not exceeded with probability at least $\alpha$. The corresponding reference problem is the chance-constrained, or equivalently, the VaR problem.
\begin{align}
    \min_{\substack{\INVESTMENT \in \setX \\ \TOC \ge 0}} \quad & \COSTX^T \INVESTMENT + \TOC \nonumber \\
    \text{s.t.} \quad
    & \mathbb{P}_{\LINESTATUS \in \setO} \Big( \Qx{\INVESTMENT}{\LINESTATUS} \le \TOC \Big) \ge \alpha. \label{problem:var}
\end{align}
For any fixed $\INVESTMENT$, the smallest feasible $\TOC$ equals $\VaR_{\alpha}(\Qx{\INVESTMENT}{\LINESTATUS})$, so~\eqref{problem:var} minimizes the investment cost plus the $\alpha$-quantile of the operational cost. Throughout this paper, $\alpha$ denotes this probability level, that is, the required probability that the operational cost under a random outage realization does not exceed $\TOC$. It describes the cost distribution and should not be confused with a statistical confidence level.

Problem~\eqref{problem:var} cannot be solved directly. $\VaR_{\alpha}$ is, in general, a nonconvex function of $\INVESTMENT$, and $\mathbb{P}$ is defined over $2^{|\setL| \times |\setT|}$ realizations, so the probability in~\eqref{problem:var} cannot be evaluated exactly. We therefore follow the scenario approach to chance-constrained optimization \citep{calafiore2005uncertain, campi2008exact}. The chance constraint is replaced by a robust constraint over a finite set $\SampledScenarios \subset \setO$ of outage realizations drawn independently from $\mathbb{P}$ (Section~\ref{sec:meth_simulation}), together with the \emph{nominal scenario} $\LINESTATUS^{\text{base}} = \bm{1}^{|\setL| \times |\setT|}$, in which all lines are always available. This results in the following optimization model. 
\begin{align}
    \min_{\substack{\INVESTMENT \in \setX \\ \TOC \ge 0}} \quad & \COSTX^T \INVESTMENT + \TOC \nonumber \\
    \text{s.t.} \quad
    & \Qx{\INVESTMENT}{\LINESTATUS} \le \TOC \qquad \forall \LINESTATUS \in \SampledScenarios \cup \{\LINESTATUS^{\text{base}}\}. \label{problem:sampled_robust}
\end{align}
In~\eqref{problem:sampled_robust}, $\TOC$ represents the worst-case operational cost across the sampled realizations. Because $\Qx{\cdot}{\LINESTATUS}$ is a convex, piecewise-linear function of $\INVESTMENT$ for every $\LINESTATUS$, \eqref{problem:sampled_robust} is a convex program for any finite $\SampledScenarios$. We call~\eqref{problem:sampled_robust} the \emph{scenario-based robust problem}. The solution method in Sections~\ref{sec:meth_simulation} and~\ref{sec:benders_decomposition} solves this problem and enlarges $\SampledScenarios$ step by step through Monte Carlo sampling. Unlike budget-based robust optimization, it requires no uncertainty set chosen in advance. Its conservatism follows from the outage process and the number of sampled realizations, and it increases as more realizations are sampled.
\subsection{Simulation of Transmission Line Outages} \label{sec:meth_simulation}
The status of each transmission line is modeled as a two-state discrete-time Markov chain with the states \emph{in service} ($1$) and \emph{out of service} ($0$). Transitions between these states depend on two reliability parameters: the base failure rate $\FailRateBase$ (failures per unit length per hour) and the mean time to repair $\MTTR$ \citep{7_Beyza_characterising}. Consistent with standard reliability models, and as described in \citet{19_yang2019failure}, we assume that the failure rate of a line is proportional to its length $\Lline{l}$.

In continuous time, line $l$ fails at the constant rate $\FailRate{l} = \FailRateBase \Lline{l}$ and is repaired at the rate $\RepairRate = 1/\MTTR$. The optimization model, however, uses the status of a line only once per time step of duration $\TIMESTEPDURATION$ (see Constraints~\eqref{eq:op_line_cap} and~\eqref{eq:op_dc_flow} in~\ref{app:subproblem_model}). We therefore use the exact per-step transition probabilities of this two-state process \citep{31_allan2013reliability, 29_Ge_2025_Continuous},
\begin{align}
    \ProbFail{l} &= \frac{\FailRate{l}}{\FailRate{l} + \RepairRate} \Big(1 - \exp\!\big(-(\FailRate{l} + \RepairRate)\,\TIMESTEPDURATION\big)\Big), \nonumber \\
    \ProbRep{l} &= \frac{\RepairRate}{\FailRate{l} + \RepairRate} \Big(1 - \exp\!\big(-(\FailRate{l} + \RepairRate)\,\TIMESTEPDURATION\big)\Big), \label{eq:transition_probabilities}
\end{align}
where $\ProbFail{l}$ denotes the probability that a line in service fails, and $\ProbRep{l}$ denotes the probability that a line out of service is repaired from one time step to the next. Each outage scenario is represented as a binary availability matrix $\LINESTATUS$ generated by Monte Carlo simulation. The status of each line in the first time step is sampled from its steady-state distribution, where the line is in service with probability $\ProbInit{l} = \RepairRate/(\FailRate{l}+\RepairRate)$. For subsequent time steps, the status of the lines is sampled using the probabilities in~\eqref{eq:transition_probabilities}. Function~\ref{func:Simulate_Failures} in~\ref{app:simulation} provides the complete procedure for generating a set of scenarios.

\subsection{Memory-Augmented Benders Decomposition} \label{sec:benders_decomposition}
Solving thousands of operational scenarios in one model is computationally intractable. We therefore use Benders decomposition, which is widely used for capacity expansion problems with recourse \citep{8_Wang2022_dist, 4_Yuan2022_TEP, 5_Moreira2015_ARO, 9_Wu2018_robust, 20_Paula2020, 18_Zhang2021_data, 21_Wen2022_robust}. The method splits the problem into a master problem (MP) and sub-problems (SPs). The MP optimizes the first-stage investments $\INVESTMENT$ subject to the technical constraints $\Ax \INVESTMENT \le \bx$. Each SP evaluates the second-stage operation $\OPERATION$ and its cost for fixed investments under one outage scenario. The SPs are independent of each other, so we solve them in parallel, which makes it possible to handle several thousand scenarios. The complete formulations of the MP and SPs are given in~\ref{app:master_model} and~\ref{app:subproblem_model}, respectively. Algorithm~\ref{alg:Benders_Decomposition} outlines the complete procedure of the proposed decomposition algorithm, and Figure~\ref{fig:benders_flowchart} in~\ref{app:benders_details} shows it as a flowchart. The procedure has two phases. Phase~I is the memory-augmented Benders decomposition, which generates candidate plans and selects an incumbent plan $\INVESTMENT^*$ (Lines~\ref{line:init_bounds}--\ref{line:loop_end}). Phase~II evaluates $\INVESTMENT^*$ by out-of-sample simulation (Lines~\ref{line:eval_scenarios}--\ref{line:return}).

\begin{algorithm*}[tp]
	\caption{BendersDecomposition($\noIterations, \noRandomScenarios, \noMemoryScenarios, \noSimScenarios, \ConfidenceLevel$)} \label{alg:Benders_Decomposition}
	\LinesNumbered
    \DontPrintSemicolon

    \tcp{Phase I: memory-augmented Benders decomposition (selection of the incumbent plan)}
    $\LB \gets -\infty, \; \UB \gets +\infty, \; \LINESTATUS^\text{base} \gets \bm{1}^{|\setL| \times |\setT|}$ \label{line:init_bounds} \tcp*{Initialize bounds and nominal scenario}
    $\MemoryScenarios \gets \text{\ref{func:Simulate_Failures}}(\noMemoryScenarios, \FailRateBase, \MTTR)$ \label{line:init_memory} \tcp*{Fill memory with an independent batch of scenarios}
        $\Pi(\LINESTATUS) \coloneqq \left\{ \DualVariables \le \bm{0} \mid \Wy(\LINESTATUS)^T \DualVariables \le \COSTy \right\}$ for $\LINESTATUS \in \setO$ \label{line:init_dual} \tcp*{Dual feasible set of the sub-problem for scenario $\LINESTATUS$}
    $\setX \gets \left\{ (\INVESTMENT, \TOC) \mid \INVESTMENT \ge \bm{0}, \Ax \INVESTMENT \le \bx, \exists \OPERATION \ge \bm{0}: \Wy(\LINESTATUS^\text{base}) \OPERATION \le \hy - \Ty{\LINESTATUS^\text{base}} \INVESTMENT, \TOC \ge \COSTy^T \OPERATION \right\}$ \label{line:init_primal} \tcp*{Initialize master problem's feasible set with nominal scenario}
    
    \For{$k \gets 1$ {\bf to} $\noIterations$ \label{line:loop_start}}{
    
        $(\overline{\INVESTMENT}, \LB) \gets \arg\min_{(\INVESTMENT, \TOC) \in \setX} \left\{\COSTX^T \INVESTMENT + \TOC \right\}$ \label{line:solve_mp} \tcp*{Solve master problem \& update $\LB$}
        
        $\RandomScenarios \gets \text{\ref{func:Simulate_Failures}}(\noRandomScenarios, \FailRateBase, \MTTR), \quad \overline{\setO} \gets \RandomScenarios \cup \MemoryScenarios$ \label{line:scenarios} \tcp*{Generate and aggregate scenarios}
        
        \tcp{Solve sub-problems and calculate system costs in parallel}
        \ParFor{$\LINESTATUS \in \overline{\setO}$ \label{line:sp_start}}{
            $\overline{\DualVariables}(\LINESTATUS) \gets \arg\max_{\DualVariables \in \Pi(\LINESTATUS)} \left\{ \DualVariables^T (\hy - \Ty{\LINESTATUS} \overline{\INVESTMENT}) \right\}$ \label{line:solve_sp}\; 
            $\TC(\LINESTATUS) \gets  \COSTX^T \overline{\INVESTMENT} + {\overline{\DualVariables}(\LINESTATUS)}^T (\hy - \Ty{\LINESTATUS} \overline{\INVESTMENT})$ \label{line:calc_tc}\;
        } \label{line:sp_end}
        
        \tcp{Memory update: keep the $\noMemoryScenarios$ most costly scenarios of $\overline{\setO}$}
        $\MemoryScenarios \gets \arg\operatorname{top}_{\noMemoryScenarios}\big\{\TC(\LINESTATUS) \mid \LINESTATUS \in \overline{\setO}\big\}$ \label{line:mem_update}\;
        
        $\setX \gets \setX \bigcap_{\LINESTATUS \in \overline{\setO}}\left\{ (\INVESTMENT, \TOC) \mid \TOC \ge {\overline{\DualVariables}(\LINESTATUS)}^T (\hy - \Ty{\LINESTATUS} \INVESTMENT)\right\}$ \label{line:add_cuts} \tcp*{Add optimality cuts}
        
        $\overline{\UB} \gets \text{Mean}_{\LINESTATUS \in \overline{\setO}}(\TC(\LINESTATUS)) + \sqrt{\frac{\ConfidenceLevel}{1- \ConfidenceLevel} \text{Var}_{\LINESTATUS \in \overline{\setO}}(\TC(\LINESTATUS))}$ \label{line:cantelli} \tcp*{In-loop Cantelli bound}

        \If{$\overline{\UB} < \UB$ \label{line:ub_update_start}}{
            $\UB \gets \overline{\UB}, \quad \INVESTMENT^* \gets \overline{\INVESTMENT}$ \label{line:ub_update} \tcp*{Update best upper bound and investment plan}
        } \label{line:ub_update_end}
        $\Gap \gets \frac{\UB - \LB}{\UB} \times 100$ \label{line:gap} \tcp*{In-loop relative gap}
    } \label{line:loop_end}
    \tcp{Phase II: out-of-sample evaluation of the incumbent plan $\INVESTMENT^*$}
    $\RandomScenarios \gets \text{\ref{func:Simulate_Failures}}(\noSimScenarios, \FailRateBase, \MTTR)$ \label{line:eval_scenarios} \tcp*{New scenarios, independent of Phase I}
    \ParFor{$\LINESTATUS \in \RandomScenarios$ \label{line:eval_start}}{
        $\TC(\LINESTATUS) \gets \COSTX^T \INVESTMENT^* + \Qx{\INVESTMENT^*}{\LINESTATUS}$ \label{line:eval_tc} \tcp*{Solve sub-problem for the fixed plan $\INVESTMENT^*$}
    } \label{line:eval_end}
    $\UB \gets \text{Mean}_{\LINESTATUS \in \RandomScenarios}(\TC(\LINESTATUS)) + \sqrt{\frac{\ConfidenceLevel}{1- \ConfidenceLevel} \text{Var}_{\LINESTATUS \in \RandomScenarios}(\TC(\LINESTATUS))}$ \label{line:ub_sim} \tcp*{Simulation-based bound}
    $\Gap \gets \frac{\UB - \LB}{\UB} \times 100$ \label{line:gap_sim} \tcp*{Reported gap}
 	\Return{$(\LB, \UB, \Gap, \INVESTMENT^*)$} \label{line:return}
\end{algorithm*}

At initialization (Lines~\ref{line:init_bounds}--\ref{line:init_primal}), the algorithm fills the scenario memory, defines the dual feasible sets $\Pi(\LINESTATUS)$, and initializes the MP feasible set $\setX$. In a standard Benders implementation, the first MP contains no information about system operation and often returns arbitrary or zero-capacity plans. We therefore add the operational constraints and the recourse cost of the nominal scenario $\LINESTATUS^\text{base}$ (Section~\ref{sec:robust_model}) directly to $\setX$ (Line~\ref{line:init_primal}). This adds little computational effort but gives a good initial plan and a tight initial lower bound.

In each iteration (Lines~\ref{line:loop_start}--\ref{line:loop_end}), the MP returns a candidate plan $\overline{\INVESTMENT}$ and a lower bound $\LB$ (Line~\ref{line:solve_mp}). A new batch of $\noRandomScenarios$ outage scenarios $\RandomScenarios$ is simulated and combined with the memory $\MemoryScenarios$ (Line~\ref{line:scenarios}). The sub-problems for all scenarios in $\overline{\setO} = \RandomScenarios \cup \MemoryScenarios$ are then solved in parallel (Lines~\ref{line:sp_start}--\ref{line:sp_end}). Involuntary loss of load is always allowed at a high penalty, so the sub-problems have complete recourse and only optimality cuts are needed (Line~\ref{line:add_cuts}; \ref{app:optimality_cut}). Because each optimality cut underestimates $\Qx{\INVESTMENT}{\LINESTATUS}$, the MP is a relaxation of~\eqref{problem:sampled_robust} for the scenarios sampled so far, and $\LB$ is a valid, nondecreasing lower bound on its optimal value.

With iterative sampling, the algorithm may \emph{forget} severe outage scenarios that are not evaluated again in later iterations. To prevent this, the algorithm keeps a memory $\MemoryScenarios$ of fixed size $\noMemoryScenarios$. The memory is initialized with an independent batch of simulated scenarios (Line~\ref{line:init_memory}). After solving the subproblems, the memory is replaced with the $\noMemoryScenarios$ scenarios in $\overline{\setO}$ with the highest total system cost $\TC(\LINESTATUS)$ (Line~\ref{line:mem_update}). The memory holds the outages that are most severe for the current plan, and each new candidate plan is tested against them. Every iteration therefore solves exactly $\noRandomScenarios + \noMemoryScenarios$ sub-problems, which keeps the parallel workload constant.

The probability level $\ConfidenceLevel$ enters only through the upper bound. A reliable upper bound is hard to obtain from samples. The empirical maximum cost depends strongly on the sample size and may underestimate the cost of unobserved realizations. An empirical $\ConfidenceLevel$-quantile, in turn, is unreliable for $\ConfidenceLevel$ close to one at moderate sample sizes. We therefore use Cantelli's inequality, the one-sided version of Chebyshev's inequality \citep{34_Barrieu_2015, 35_Mercadier_2021, 36_apollonio2025cantelli}. For a candidate plan $\overline{\INVESTMENT}$, let $\mu_{\TC}$ and $\sigma_{\TC}^2$ denote the mean and variance of the total system cost $\TC(\LINESTATUS)$ under $\mathbb{P}$. Cantelli's inequality states that $\mathbb{P}\big(\TC(\LINESTATUS) \ge \mu_{\TC} + \kappa\big) \le \sigma_{\TC}^2/(\sigma_{\TC}^2 + \kappa^2)$ for any $\kappa > 0$. Choosing $\kappa = \sigma_{\TC}\sqrt{\ConfidenceLevel/(1-\ConfidenceLevel)}$ gives
\begin{multline}
    \mathbb{P}\Big( \TC(\LINESTATUS) < \mu_{\TC} + \sqrt{\tfrac{\ConfidenceLevel}{1-\ConfidenceLevel}}\, \sigma_{\TC} \Big) \ge \ConfidenceLevel \\
    \text{and hence} \quad
    \VaR_{\ConfidenceLevel}\big(\TC(\LINESTATUS)\big) \le \mu_{\TC} + \sqrt{\tfrac{\ConfidenceLevel}{1-\ConfidenceLevel}}\, \sigma_{\TC}. \label{eq:cantelli}
\end{multline}

The total system cost $\TC(\LINESTATUS)$ differs from $\Qx{\overline{\INVESTMENT}}{\LINESTATUS}$ only by the constant investment cost $\COSTX^T \overline{\INVESTMENT}$. The right-hand side of~\eqref{eq:cantelli} is therefore an upper bound on the objective value of $\overline{\INVESTMENT}$ in the VaR problem~\eqref{problem:var}, and thus also on the optimal value of~\eqref{problem:var}. In Line~\ref{line:cantelli}, $\mu_{\TC}$ and $\sigma_{\TC}^2$ are replaced by their sample estimates over $\overline{\setO}$, which gives the in-loop bound $\overline{\UB}$. The incumbent plan $\INVESTMENT^*$ is the candidate with the lowest $\overline{\UB}$ found so far (Lines~\ref{line:ub_update_start}--\ref{line:ub_update_end}).

The in-loop bound $\overline{\UB}$ is estimated from only $\noRandomScenarios + \noMemoryScenarios$ scenarios per iteration, and the memory scenarios among them are selected for their high cost. It is therefore a biased, approximate bound and not a formal statistical guarantee (\ref{app:benders_details}). In Phase~I, we use it only to select the incumbent plan and to monitor convergence through the in-loop gap (Line~\ref{line:gap}).

After $\noIterations$ iterations, Phase~II evaluates $\INVESTMENT^*$ on a new and much larger set $\RandomScenarios$ of $\noSimScenarios$ outage scenarios, generated independently of all scenarios used in Phase~I (Line~\ref{line:eval_scenarios}). With $\INVESTMENT^*$ fixed, the sub-problems are solved in parallel (Lines~\ref{line:eval_start}--\ref{line:eval_end}), and the resulting system costs give the reported $\UB$ and $\Gap$ (Lines~\ref{line:ub_sim}--\ref{line:return}). These scenarios contain no memory scenarios and play no role in selecting $\INVESTMENT^*$. Their sample moments are therefore unbiased estimates for $\INVESTMENT^*$, as in standard Monte Carlo procedures for assessing solution quality \citep{mak1999monte}. This removes the upward bias caused by the memory scenarios and the downward bias caused by selecting the candidate with the lowest of many noisy estimates \citep{smith2006optimizer}. The larger sample also reduces the sampling error. The reported $\UB$ nevertheless remains approximate, because it is still computed from estimated moments. In addition, $\LB$ and $\UB$ refer to two related but different quantities. The lower bound $\LB$ applies to the worst-case cost over the sampled scenarios that any plan can attain, while the upper bound $\UB$ applies to the $\VaR_{\ConfidenceLevel}$ of the incumbent plan's cost. The $\Gap$ is therefore not a formal optimality certificate for either~\eqref{problem:var} or~\eqref{problem:sampled_robust}, but an indicator of solution quality.

\section{Case Study} \label{sec:case_study}

This section describes the case study. Section~\ref{sec:case_system} presents the test system and its data, Section~\ref{sec:outage_cases} defines the outage cases, and Section~\ref{sec:benchmarks} introduces the benchmarks and the metrics used to evaluate the capacity plans.

\subsection{Test System and Data} \label{sec:case_system}

We apply the proposed framework to a reduced European power system with techno-economic parameters for 2030 and weather and demand data from 2013. The network consists of 25 nodes covering Germany and its interconnected neighboring countries (Austria, Belgium, Switzerland, the Czech Republic, Denmark, France, Luxembourg, the Netherlands, Poland, and Sweden). Since Germany is the focus of this study, we give it a higher spatial resolution than its neighboring countries. To do so, we assign a high clustering weight to Germany and low clustering weights to the neighboring countries in \texttt{PyPSA-Eur} (Table~\ref{tab:sys_params}). Because the clustering keeps at least one node per country and synchronous area, each neighboring country is represented by one node and Germany by 14 nodes. The only exception is Denmark with two nodes, as western and eastern Denmark belong to different synchronous areas (Continental Europe and the Nordic system).

Generation and storage capacities are fully greenfield and can be built from zero at every node, as in the formulation in Section~\ref{sec:ProblemStatement}. The transmission network, in contrast, builds on the existing European grid. For each of the 56 lines of a 25-node network generated with the open-source \texttt{PyPSA-Eur} model \citep{horsch2018pypsa} (49 AC lines and 7 DC links; see~\ref{app:data}), we set its existing capacity as the lower investment bound $\mincaptrans{l}$ in Constraint~\eqref{eq:inv_trans_potential}. The upper bound $\maxcaptrans{l}$ is the existing capacity plus the default expansion limit of \texttt{PyPSA-Eur} (\SI{20}{\GW} for AC lines and \SI{30}{\GW} for DC links). The annualized investment costs $\ANNUALCOSTLine{l}$ are also taken from \texttt{PyPSA-Eur}. As stated in Assumption~\ref{ass:corridor}, each line may aggregate several parallel circuits, and failure statistics are applied at the corridor level. An outage in the model therefore means that the full capacity of a line is lost at once.

The technology portfolio comprises three renewable sources (solar, onshore wind, and offshore wind), three dispatchable sources (nuclear, combined-cycle gas turbines, and biogas), and hydrogen storage, with techno-economic parameters also projected for 2030. To restrict emissions, we set a carbon cap of \SI{0.005}{\TonCOtwo \per \MWh}. For computational tractability, the 8760 hourly load and renewable capacity factor profiles are downsampled to 25-hour time steps. Each step takes the value of the first hour of its block, which keeps the hourly variability of the original profiles instead of averaging it out. Because each step is one hour longer than a day, the sampled hours also move through all hours of the day. This gives $|\setT| = 350$ time steps, and operational costs are scaled to a full year with the weight $\OPWEIGHT = 8760/350$ (\ref{app:subproblem_model}). At this resolution, the model cannot represent intraday charge--discharge cycles, which provide most of the value of short-duration batteries. We therefore exclude battery storage and model only long-duration hydrogen storage. Hydrogen storage typically charges and discharges over several days or weeks, so 25-hour time steps resolve its operation adequately. All techno-economic parameters are given in~\ref{app:data}.

\subsection{Outage Cases} \label{sec:outage_cases}

The simulation of line failure and repair uses two standard reliability parameters (Section~\ref{sec:meth_simulation}). The first is the line failure rate $\FailRateBase$ in failures per $100$~km per year, which is converted to failures per km per hour with $\HOURSPERYEAR = 8760$~h. The second is the mean time to repair (MTTR). Empirical studies report outage rates of about $0.5$ to $2$ failures per $100$~km per year for overhead high-voltage transmission lines, depending on network characteristics, environmental conditions, and voltage level \citep{30_pijarski2025outage}. Reported repair times also vary across systems and commonly range from several hours to about one day \citep{30_pijarski2025outage,33_entsoe2024nordic}. Based on these ranges, we define three cases with increasing outage severity, summarized in Table~\ref{tab:test_cases}.

\begin{table}[htbp]
    \centering
    \begin{threeparttable}
        \caption{Definition of the outage cases.}
        \label{tab:test_cases}
        \small\setlength{\tabcolsep}{4pt}
        \begin{tabular}{l l c c}
            \toprule
            \textbf{Outage Case} & \textbf{Description} & \textbf{MTTR} & \textbf{Failure Rate} \\
            & & \textbf{(h)} & \textbf{(per 100 km/yr)} \\
            \midrule
            \textbf{Case A} & Low Severity    & 12 & 0.5 \\
            \textbf{Case B} & Medium Severity & 24 & 1.0 \\
            \textbf{Case C} & High Severity   & 48 & 2.0 \\
            \bottomrule
        \end{tabular}
    \end{threeparttable}
\end{table}

\subsection{Benchmarks and Evaluation Metrics} \label{sec:benchmarks}

We evaluate the proposed model by comparing its capacity plans with those of two benchmarks, its deterministic counterpart and a reserve-margin approach, which is common practice in energy system modeling. All three capacity plans are then stress-tested on the same simulated outage scenarios.

The deterministic benchmark is the G\&TEP model in~\ref{app:master_model} that considers only the nominal scenario $\LINESTATUS^{\text{base}}$ (Section~\ref{sec:robust_model}), in which all transmission lines are always available. It therefore finds an optimal capacity mix under the $N$-0 security assumption. The reserve-margin heuristic is a well-known approach for handling uncertainty in energy system planning, in which an exogenous reserve margin is imposed \citep{stephen2022enhanced, 46_giglio_2025_integrated}. It is the deterministic model with the reserve margin constraint activated, which requires the total firm capacity at \emph{each node} to be at least $\RESERVEMARGIN$ times its peak demand (Constraint~\eqref{eq:inv_res_margin}). Following \citet{47_pypsa_user_meeting_2024} and \citet{46_giglio_2025_integrated}, we use $\RESERVEMARGIN = 1.1$. The capacity credit $\CAPACITYCREDIT{g}$ of each technology in this constraint is given in Table~\ref{tab:tech_params}.\footnote{With $\RESERVEMARGIN \ge 1$, the demand at each node can be met by local firm capacity alone, as if the node were cut off from the network. This reduces LoL when transmission lines fail.}

All capacity plans are evaluated as in Phase~II of Algorithm~\ref{alg:Benders_Decomposition}, that is, on $\noSimScenarios = 1000$ simulated outage scenarios per outage case, which are independent of the scenarios used during optimization. We assess cost performance by the highest, lowest, and mean system costs in these simulations. Robustness is assessed mainly by the loss of load (LoL) and also by the standard deviation of the system cost.

In the optimization, LoL carries a high penalty of \SI{100}{\Currency\per\kWh}. This penalty does not represent the actual economic cost of unserved energy. Its purpose is to make the model meet demand whenever possible instead of using LoL as an easy way to stay feasible. Because the actual cost of unserved energy is hard to estimate, we do not use it in the reported cost. Therefore, for a fair comparison that is not distorted by the penalty, we report an \emph{adjusted system cost}~\eqref{eq:adjusted_cost}, in which LoL is priced at the same rate as voluntary load shedding, $\COSTSHED = \SI{1}{\Currency\per\kWh}$.
\begin{multline}
    \text{Adjusted System Cost} = \text{Total System Cost} \\
    - \OPWEIGHT (\COSTLOSS - \COSTSHED)  \sum_{n \in \setN}\sum_{t \in \setT} \Ploss{n}{t}. \label{eq:adjusted_cost}
\end{multline}

We report LoL separately in \si{\GWh} so that the cost and robustness of each plan can be assessed independently. We measure LoL in two ways. The first is how often it occurs, given by the share of scenarios with any LoL. Because each scenario covers one year, this share is the probability of at least one LoL event in a year. The second is how much load is lost, given by the mean LoL over all scenarios and the highest LoL in any scenario.

The decomposition algorithm is run for $\noIterations = 300$ iterations, with $\noRandomScenarios = 100$ new random scenarios per iteration and $\noMemoryScenarios = 30$ scenarios in memory. All computations are performed on a workstation with an Intel Core i9-14900 CPU (24 cores, 32 threads), 64~GB RAM, and Windows 11 Enterprise. The models are implemented in Julia and solved with Gurobi Optimizer~\citep{gurobi} version~12.0.0.

\section{Results} \label{sec:results}

This section presents the results of the case study. Section~\ref{sec:methods} addressed RQ\ref{rq:integration} by formulating a scenario-based robust optimization model that includes transmission line outages, characterized by their rate and duration, in an inter-regional G\&TEP. Sections~\ref{sec:deterministic_vs_robust}--\ref{sec:sensitivity_fr_mttr} address RQ\ref{rq:comparison}. Section~\ref{sec:deterministic_vs_robust} compares the cost and LoL of the deterministic, reserve-margin, and robust plans, Section~\ref{sec:cost_trade_off} examines the trade-off between cost and robustness and how the robust plans adapt the capacity mix, and Section~\ref{sec:sensitivity_fr_mttr} analyzes how the failure rate and the repair time affect LoL. Section~\ref{sec:algorithm_performance} addresses RQ\ref{rq:solvability} by evaluating the computational performance of the algorithm. The detailed numerical results are given in~\ref{app:detailed_results}. Table~\ref{tab:results_investment} reports the investment decisions of each plan, and Table~\ref{tab:results_cases} reports the bounds, computation times, and simulated cost and LoL statistics.

\subsection{Comparison of the Deterministic, Reserve-Margin, and Robust Plans} \label{sec:deterministic_vs_robust}

\begin{figure*}[tp]
    \centering
    \includegraphics[width=0.8\textwidth]{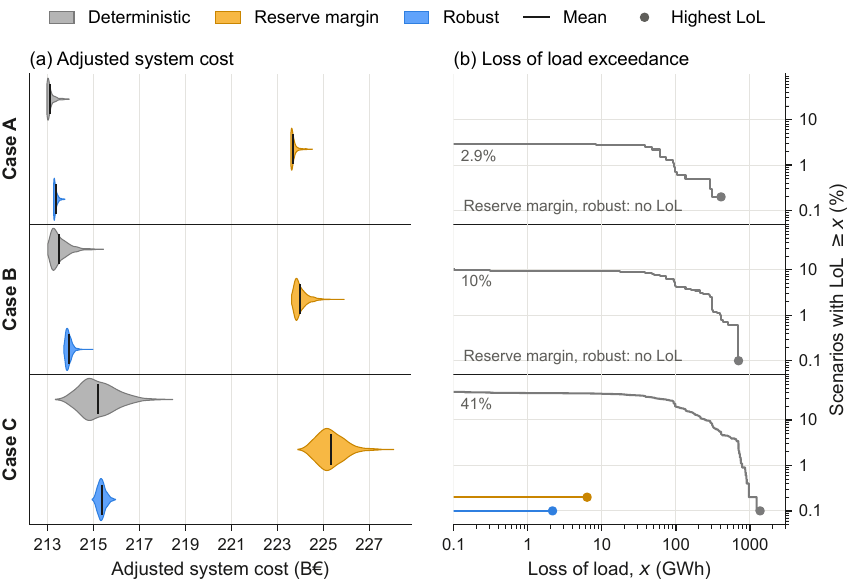}
    \caption{Adjusted system cost (a) and LoL (b) of the deterministic, reserve-margin, and robust plans in outage Cases~A, B, and~C, based on 1000 out-of-sample scenarios. In (a), violins show the distribution of the adjusted system cost and vertical lines show the mean. In (b), each curve shows the share of scenarios in which LoL is at least the value on the horizontal axis. The left end of a curve is the share of scenarios with any LoL, and the dot marks the highest LoL.}
    \label{fig:violin_cost_lol}
\end{figure*}

We first compare the three capacity plans by stress-testing them on the same 1000 out-of-sample outage scenarios in each outage case (Section~\ref{sec:benchmarks}). Figure~\ref{fig:violin_cost_lol} shows the distributions of their adjusted system cost (a) and LoL (b), and Table~\ref{tab:results_cases} gives the corresponding statistics.

As seen in Figure~\ref{fig:violin_cost_lol}b, the deterministic plan is vulnerable to transmission outages. Without enough local backup or transmission redundancy, line outages restrict power flows and lead to LoL. LoL becomes both more frequent and larger as outages become more severe. The share of scenarios with LoL, that is, the probability of at least one LoL event in a year, rises from about 3\% in Case~A to 41\% in Case~C. In Case~C, the most severe case, the mean LoL over all scenarios is about 79~\si{\GWh}, and the highest LoL is 1371~\si{\GWh} (about seven hours of average system load). The deterministic plan also has the largest cost variability of the three approaches in all outage cases, as seen in Figure~\ref{fig:violin_cost_lol}a.

The reserve-margin heuristic reduces this vulnerability by building local firm generation. It almost eliminates LoL; as seen in Figure~\ref{fig:violin_cost_lol}b, no LoL occurs in Cases~A and~B, and the highest LoL in Case~C is 6~\si{\GWh}. However, this reliability is costly. In Case~C, the mean adjusted system cost of the heuristic is 10.1~\si{B\Currency} (4.7\%) higher than that of the deterministic plan (Figure~\ref{fig:violin_cost_lol}a). The main reason is the investment in redundant local capacity, which raises the annualized investment and fixed O\&M cost by 8.5~\si{B\Currency}. The mean operational cost is also 2.7~\si{B\Currency} higher (Table~\ref{tab:results_cases}). Both differences follow from the capacity mix of the heuristic, which we examine in Section~\ref{sec:cost_trade_off}.

The robust framework, in contrast, accounts for line outages directly in the investment decisions. The results from 1000 simulated scenarios show that it can also almost eliminate LoL. No LoL occurs in Cases~A and~B, and the highest LoL in Case~C is just 2~\si{\GWh} (Figure~\ref{fig:violin_cost_lol}b). In all three cases, its mean adjusted system cost is 4.4--4.6\% lower than that of the heuristic and only 0.1--0.2\% higher than that of the deterministic plan (Figure~\ref{fig:violin_cost_lol}a). It should be mentioned that this comparison is conservative because the adjusted system cost values LoL at only \SI{1}{\Currency\per\kWh}. The robust plan would have a lower mean cost than the deterministic plan if lost load were valued at more than about \SI{3}{\Currency\per\kWh}, \SI{24}{\Currency\per\kWh}, and \SI{70}{\Currency\per\kWh} in Cases~C, B, and~A, respectively.

Figure~\ref{fig:violin_cost_lol} also shows that the robust plan has the narrowest cost distribution of the three approaches in all cases. In Case~C, the standard deviation of its adjusted system cost is 0.2~\si{B\Currency}, compared with 0.6~\si{B\Currency} for the heuristic and 0.8~\si{B\Currency} for the deterministic plan. Its highest adjusted system cost is also 2.5~\si{B\Currency} lower than that of the deterministic plan (see Table~\ref{tab:results_cases}).

\subsection{Cost-Robustness Trade-Offs and Capacity Adaptations} \label{sec:cost_trade_off}

\begin{figure*}[tp]
    \centering
    \includegraphics[width=0.8\textwidth]{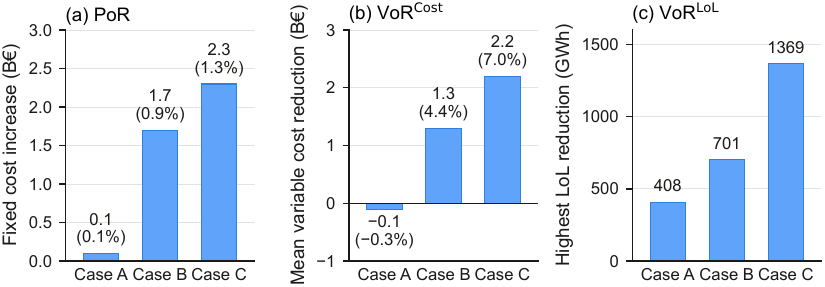} %
    \caption{Price and value of robustness in outage Cases~A, B, and~C, based on 1000 out-of-sample scenarios per case. (a) PoR, the increase in annualized investment and fixed O\&M cost. (b) \VoRCost, the reduction in mean variable cost, that is, the adjusted system cost without annualized investment and fixed O\&M costs. The negative value in Case~A denotes a slight increase. (c) \VoRLoL, the reduction in the highest LoL across the simulated scenarios. Percentages are relative to the deterministic plan.}
    \label{fig:PoR_VoR}
\end{figure*}

To quantify the trade-off between cost and robustness, we use two common metrics, the \emph{price of robustness} (PoR) and the \emph{value of robustness} (VoR). The PoR is the additional fixed investment that the robust model needs, relative to the deterministic plan, to protect against transmission outages. In the mild Case~A, the PoR is only 0.1~\si{B\Currency}, which corresponds to a 0.1\% increase over the deterministic plan (Figure~\ref{fig:PoR_VoR}a). In the more severe Cases~B and~C, it rises to 1.7~\si{B\Currency} and 2.3~\si{B\Currency}, by respectively 0.9\% and 1.3\%.

We measure the benefits of this additional investment with two VoR metrics. The first, \VoRCost, is the reduction in mean variable cost relative to the deterministic plan, where the variable cost is the adjusted system cost without the annualized investment and fixed O\&M costs. The second metric, \VoRLoL, represents the reduction in the highest Loss of Load (LoL) across the simulated scenarios. In Case~A, the robust plan has a slightly higher variable cost than the deterministic plan, resulting in a slightly negative \VoRCost (Figure~\ref{fig:PoR_VoR}b). In contrast, in Cases~B and~C, the robust plan reduces the variable cost by 4.4\% and 7.0\%, respectively. These reductions offset most of the additional investment required. The cost savings arise from two primary sources (see Table~\ref{tab:results_cases}). First, operational costs decrease because the additional nuclear capacity in the robust plans has lower variable costs than the gas and biogas units it partially replaces. Second, load shedding during outages is reduced. Regarding robustness, the robust plans decrease the highest LoL by 408~\si{\GWh}, 701~\si{\GWh}, and 1369~\si{\GWh} in Cases~A, B, and~C, respectively (Figure~\ref{fig:PoR_VoR}c).

\begin{figure*}[tp]
    \centering
    \includegraphics[width=0.8\textwidth]{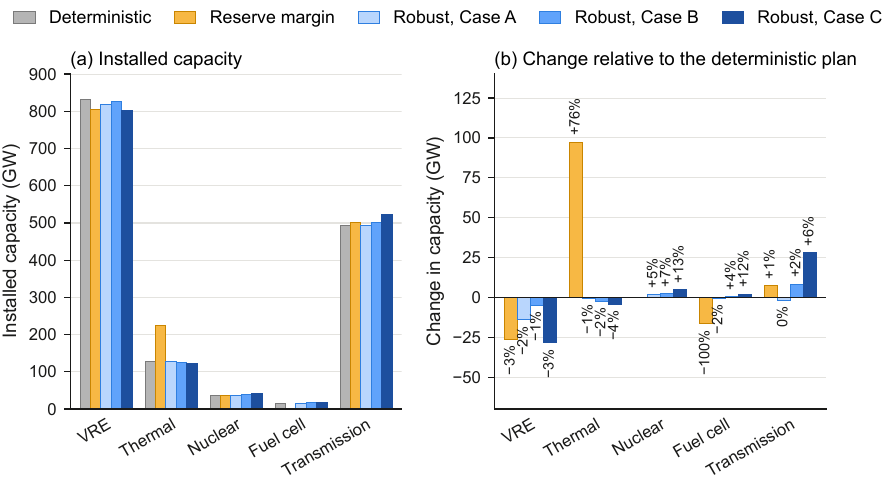} %
    \caption{Installed capacity (a) and change relative to the deterministic plan (b) of the deterministic, reserve-margin, and robust plans. The robust plans are optimized separately for outage Cases~A, B, and~C, where severity increases as both the failure rate and repair time double from one case to the next. VRE includes onshore wind, offshore wind, and solar capacity, and thermal capacity includes gas and biogas. Labels in (b) indicate the change relative to the deterministic plan in percent.}
    \label{fig:capacity_mix}
\end{figure*}

To understand how the system adapts, we next examine the capacity mix (Figure~\ref{fig:capacity_mix}a and Table~\ref{tab:results_investment}), its change relative to the deterministic plan (Figure~\ref{fig:capacity_mix}b), and its spatial distribution (Figure~\ref{fig:all_plans}). The reserve-margin heuristic protects against line outages mainly by building much more firm generation. Its thermal capacity (gas plus biogas capacity) is 76\% higher than in the deterministic plan, and it builds no hydrogen storage. Without storage to shift surplus renewable electricity in time, the heuristic covers more demand with gas and biogas. These technologies have much higher variable costs than wind and solar (Table~\ref{tab:tech_params}). This explains its higher operational cost mentioned in Section~\ref{sec:deterministic_vs_robust}. The robust plan, by contrast, slightly reduces thermal capacity by 3.5\% in Case~C. It relies more on nuclear capacity and fuel cells as outages become more severe, with 13\% and 12\% more capacity, respectively, in Case~C.

All plans expand the existing transmission grid of \SI{319}{\GW} substantially, by 54\% to 64\%. In comparison to the deterministic plan, total transmission capacity remains nearly unchanged in Case~A and increases only moderately in Cases~B and~C, by 2\% and 6\%, respectively. However, changes on individual lines are more pronounced. In Case~C, the robust plan reduces the capacity of lines connecting Sweden and Germany via Denmark by up to 13~\si{\GW}. Conversely, it increases the capacity of direct DC links from Sweden to Germany and Poland, as well as the lines between Germany and France (see Figure~\ref{fig:all_plans}, transmission row).

The robust plans also alter the geographic distribution of capacity. They reduce reliance on onshore wind, particularly in Case~C, where onshore wind capacity decreases by 7\% (see Table~\ref{tab:results_investment}). This reduction is concentrated in Sweden and is partially compensated by increased onshore wind deployment in Germany, which is closer to the main load centers (see Figure~\ref{fig:all_plans}, VRE row). One plausible explanation is that failure rates are proportional to transmission line length, making remote generation that depends on long transmission lines especially vulnerable to outages.

Hydrogen storage also changes mainly in where it is built rather than in how much. In the deterministic plan, it is installed at only five nodes (Figure~\ref{fig:all_plans}, fuel-cell row), mainly in France, Sweden, and Poland, where it shifts surplus wind generation over time. In the robust plans, fuel cells are built at 21, 24, and all 25 nodes in Cases~A, B, and~C, respectively. The share of fuel-cell capacity located in Germany rises from 4\% to 47\% in Case~C. The robust plans thus increasingly use hydrogen storage as a distributed local buffer that can supply individual nodes during transmission line outages. Conversely, gas capacity becomes concentrated at fewer nodes, with no gas capacity in Germany in Case~C.

\subsection{Sensitivity to Failure Rate and Repair Time} \label{sec:sensitivity_fr_mttr}

\begin{figure*}[tp]
    \centering
    \includegraphics[width=0.8\textwidth]{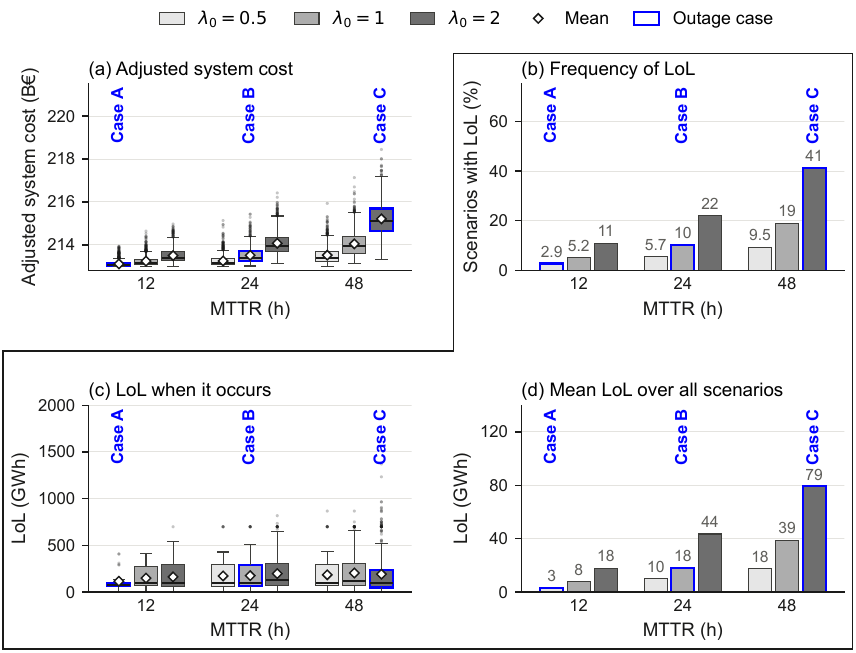}
    \caption{Sensitivity of the deterministic capacity plan to the mean time to repair (MTTR, horizontal axis) and the failure rate $\FailRateBase$ (gray shades), based on 1000 simulated outage scenarios for each combination. (a) Adjusted system cost. (b) Share of scenarios with LoL. (c) LoL in these scenarios only. (d) Mean LoL over all scenarios. Boxes span the interquartile range, horizontal lines mark the median, diamonds mark the mean, whiskers extend to 1.5 times the interquartile range, and dots are values beyond the whiskers. The framed panels (b)--(d) show LoL, and blue outlines and labels mark outage Cases~A, B, and~C.}
    \label{fig:sensitivity_fr_mttr}
\end{figure*}

In Cases~A, B, and~C, the failure rate and the mean time to repair (MTTR) increase simultaneously (Table~\ref{tab:test_cases}). The case comparison therefore cannot show which of the two parameters makes a plan more vulnerable. To separate their effects, we evaluate the deterministic plan for all nine combinations of three failure rates (0.5, 1, and 2 failures per \SI{100}{\km} per year) and three MTTRs (\SI{12}{\hour}, \SI{24}{\hour}, and \SI{48}{\hour}), with 1000 simulated outage scenarios for each combination. In terms of (MTTR, failure rate), Cases~A, B, and~C are the combinations (\SI{12}{\hour}, 0.5), (\SI{24}{\hour}, 1), and (\SI{48}{\hour}, 2), respectively. We use the deterministic plan because, unlike robust plans, its capacities do not depend on the outage parameters, so differences between combinations come solely from the outage process. Figure~\ref{fig:sensitivity_fr_mttr}a shows the adjusted system cost. Because most scenarios have no LoL, LoL is shown in three parts. Figure~\ref{fig:sensitivity_fr_mttr}b shows the share of scenarios with LoL; Figure~\ref{fig:sensitivity_fr_mttr}c shows the LoL in these scenarios only (hereafter called positive LoL); and Figure~\ref{fig:sensitivity_fr_mttr}d shows the mean LoL over all scenarios.

The results indicate that higher failure rates and longer repair times increase both the frequency and magnitude of LoL, as expected. The share of scenarios with LoL rises from about 3\% to 41\% and roughly doubles whenever either parameter doubles (Figure~\ref{fig:sensitivity_fr_mttr}b). At each failure rate, increasing MTTR from \SI{12}{\hour} to \SI{48}{\hour} raises the mean positive LoL by 18--60\% and approximately doubles the largest LoL (Figure~\ref{fig:sensitivity_fr_mttr}c). In Case~C (\SI{48}{\hour}, 2), the mean positive LoL is slightly lower than with a shorter MTTR or a lower failure rate, although the largest LoL is the highest observed (1371~\si{\GWh}). Failure frequency still has the larger effect, so the mean LoL rises from 3~\si{\GWh} to 79~\si{\GWh} mostly because LoL becomes more frequent (see Figure~\ref{fig:sensitivity_fr_mttr}b and Figure~\ref{fig:sensitivity_fr_mttr}d). The two parameters also reinforce each other. Raising both increases the mean LoL by about 2.6 times the sum of the increases from raising each alone.

Combinations with the same expected line unavailability, i.e., the same product of failure rate and MTTR, lead to similar mean outcomes but different distributions. For example, the combinations (\SI{12}{\hour}, 2), (\SI{24}{\hour}, 1), and (\SI{48}{\hour}, 0.5) yield almost the same mean LoL (about 18~\si{\GWh}, Figure~\ref{fig:sensitivity_fr_mttr}d) and mean adjusted system cost (Figure~\ref{fig:sensitivity_fr_mttr}a). As outages become rarer but longer, however, LoL occurs slightly less often (Figure~\ref{fig:sensitivity_fr_mttr}b) but is larger when it occurs, and the upper tail becomes heavier. The worst-case LoL increases from about 700~\si{\GWh} to 867~\si{\GWh} (see Figure~\ref{fig:sensitivity_fr_mttr}c, for those three combinations). The mean positive LoL likewise increases with the MTTR for the other combinations with equal unavailability.

The adjusted system cost follows the same pattern on a much smaller scale (Figure~\ref{fig:sensitivity_fr_mttr}a). Raising either parameter alone by a factor of four increases the mean cost by about 0.2\%, raising both increases it by 1.0\%, and the standard deviation of the system cost grows from 0.2~\si{B\Currency} to 0.8~\si{B\Currency}. These results show that neither parameter dominates the vulnerability of the deterministic plan in this case study. The expected line unavailability (failure rate) largely determines how often LoL occurs and its mean (Figures~\ref{fig:sensitivity_fr_mttr}b and~\ref{fig:sensitivity_fr_mttr}d), whereas longer repair times make individual LoL events larger and the upper tail heavier (Figure~\ref{fig:sensitivity_fr_mttr}c).

\subsection{Computational Performance and Convergence of the Decomposition Algorithm} \label{sec:algorithm_performance}

\begin{figure*}[tp]
    \centering
    \includegraphics[width=0.8\textwidth]{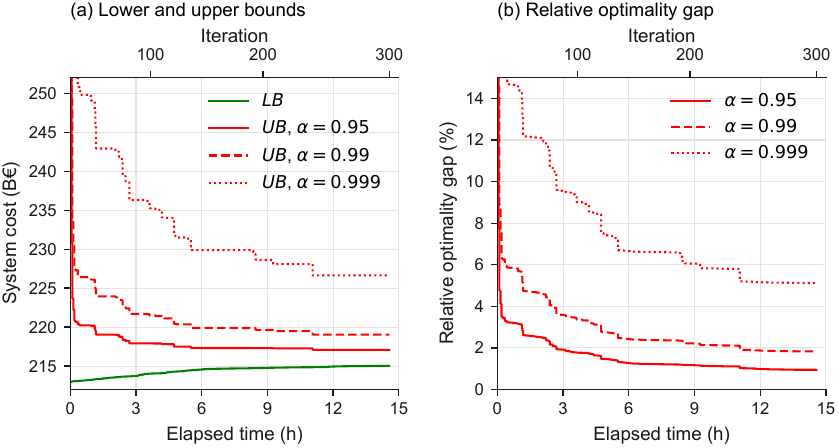}
    \caption{Convergence of Phase~I of the proposed algorithm for Case~C and three VaR probability levels, $\alpha \in \{0.95, 0.99, 0.999\}$. (a) Lower bound ($\LB$) from the master problem and in-loop upper bound ($\UB$), i.e., the Cantelli bound on the VaR of the system cost of the incumbent plan before the out-of-sample evaluation of Phase~II. (b) Relative gap $(\UB - \LB)/\UB \times 100$. The top axis gives the iteration number.}
    \label{fig:convergence}
\end{figure*}

Figure~\ref{fig:convergence}a shows the evolution of the lower bound ($\LB$) and the in-loop upper bound ($\UB$), Line~\ref{line:cantelli} of Algorithm~\ref{alg:Benders_Decomposition}, during Phase~I for the most severe Case~C. In the early iterations, the $\UB$ decreases as the algorithm identifies critical outage scenarios, adds the corresponding Benders cuts, and keeps those scenarios in memory. Meanwhile, the $\LB$ increases from about 213.0~\si{B\Currency} to 215.0~\si{B\Currency}, as the investment decisions are adjusted to reduce LoL.

Because $\alpha$ does not enter the Benders cuts, the $\LB$ trajectory is the same for all three values of $\alpha$, and only the $\UB$ differs. A higher $\alpha$ increases the Cantelli multiplier $\sqrt{\alpha/(1-\alpha)}$ (about 4.4, 9.9, and 31.6 for $\alpha = 0.95$, $0.99$, and $0.999$) and thus the $\UB$ and the gap. After 300 iterations, the in-loop gap becomes 0.9\%, 1.8\%, and 5.1\%, respectively (see Figure~\ref{fig:convergence}b).

\begin{figure*}[tp]
    \centering
    \includegraphics[width=0.8\textwidth]{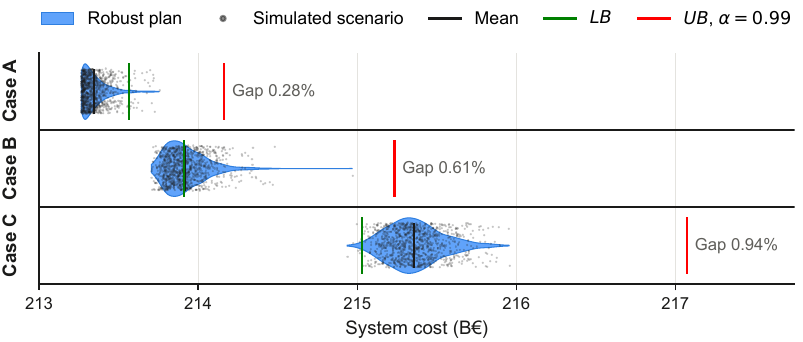} %
    \caption{Distribution of the system cost of the robust plans across 1000 out-of-sample simulations per case (Phase~II). $\LB$ is the master-problem lower bound at the end of Phase~I, and $\UB$ is the simulation-based Cantelli bound on the VaR of the system cost at $\alpha = 0.99$. The relative gap is $(\UB - \LB)/\UB \times 100$. Dots show the individual simulated costs, and black lines mark the mean.}
    \label{fig:violin}
\end{figure*}

In Phase~II (after performing 1000 simulations), the selected plans yield simulation-based upper bounds of 214.2~\si{B\Currency}, 215.2~\si{B\Currency}, and 217.1~\si{B\Currency} for Cases~A, B, and C, respectively, at $\alpha = 0.99$ (Figure~\ref{fig:violin}). Given the limitations discussed in Section~\ref{sec:benders_decomposition} and~\ref{app:benders_details}, and ignoring the sampling error in the estimated mean and standard deviation, the probability that the system cost of a selected plan exceeds these values under a random outage realization is therefore at most 1\%. For Case~C, the simulation-based $\UB$ (Line~\ref{line:ub_sim} in Algorithm~\ref{alg:Benders_Decomposition}) is 1.9~\si{B\Currency} lower than the in-loop $\UB$ (Line~\ref{line:ub_update} in Algorithm~\ref{alg:Benders_Decomposition}), which is 219.0~\si{B\Currency} (see $\UB$ in Figure~\ref{fig:convergence}a at $\ConfidenceLevel = 0.99$). This result aligns with the upward bias introduced by the memory scenarios in the in-loop bound and indicates that this bias outweighs the optimistic effect of selecting the incumbent.

With $\LB$s of 213.6~\si{B\Currency}, 213.9~\si{B\Currency}, and 215.0~\si{B\Currency}, the relative gaps (Line~\ref{line:gap_sim} of Algorithm~\ref{alg:Benders_Decomposition}) are 0.28\%, 0.61\%, and 0.94\% for Cases A, B, and C, respectively (see also Figure~\ref{fig:violin}). For Case~C, the simulation-based gap is about half the in-loop gap. The larger gaps in the more severe cases reflect greater cost variability caused by more frequent and longer outages, which widens the Cantelli bound. Although these are not formal optimality certificates (Section~\ref{sec:benders_decomposition}), the resulting gaps suggest that the conservatively bounded value at risk (VaR) of the selected plans closely approximates a lower bound on the worst-case cost achievable by any plan under the simulated scenarios. Note that the $\LB$ does not bound individual simulated costs, as it represents the worst-case scenario across all sampled cases during optimization, while each simulated value corresponds to a single scenario. Consequently, most simulated costs fall below the $\LB$ in Case~A and above it in Case~C, where the gap is largest; see Figure~\ref{fig:violin}, which shows the gaps widening as outage severity increases from Case~A to Case~C.

\section{Discussion} \label{sec:discussion}
This section discusses what the results of Section~\ref{sec:results} suggest for energy system modeling, followed by the main limitations of our approach and directions for future work.

\subsection{Implications for Energy System Modelers and Comparison with Prior Work} \label{sec:implications}

We derive three main findings from our results, subject to the limitations discussed in Section~\ref{sec:limitations}. First, treating transmission outages as static, time-independent events, an assumption implicit in $N$-0 or $N$-1 security criteria, appears to understate reliability risk, and to a lesser extent cost, when renewable variability is considered. This occurs because the operational impact of an outage depends on its timing and duration, not just on its occurrence. The sensitivity analysis in Section~\ref{sec:sensitivity_fr_mttr} supports this. Outage processes with the same expected line unavailability resulted in similar mean loss of load but different worst-case loss of load. Less frequent but longer outages produce fewer, larger events and a heavier upper tail. Therefore, outage frequency and duration should be represented separately rather than combined into a single aggregate availability measure. This limitation is common to all TEP and G\&TEP studies summarized in Table~\ref{tab:lit_review} \citep{5_Moreira2015_ARO, 9_Wu2018_robust, 18_Zhang2021_data, 8_Wang2022_dist, 4_Yuan2022_TEP}, which represent contingencies as static snapshots instead of time-dependent events. Studies that explicitly model outage duration have thus far been limited to small, single-system scales, such as microgrid scheduling under islanding \citep{24_yuanxiong_islanding, 27_Haifeng_bilevel} and sequential reliability simulation \citep{7_Beyza_characterising}. Our results presented in Section~\ref{sec:cost_trade_off} are consistent with the qualitative findings of that literature, which show that prolonged, chronologically resolved outages require local buffering resources. However, this work extends the analysis to a long-term, inter-regional planning context, which previous smaller-scale studies did not address.

Second, a comparison between the reserve-margin heuristic and the robust framework shows that rule-based heuristics, although straightforward to implement and common in the literature \citep{stephen2022enhanced, 46_giglio_2025_integrated}, can lead to inefficient technology choices. This is because such heuristics enforce uniform, node-level requirements instead of allowing the optimization process to determine where additional protection is most economically justified. Both the heuristic and robust plans nearly eliminate loss of load across all outage scenarios (with maximum losses of 6~\si{\GWh} and 2~\si{\GWh}, respectively, in the most severe simulated scenario of Case~C). The robust plans provide this level of protection at a 4.4--4.6\% (approximately 10~\si{B\Currency} per year) lower mean adjusted system cost, while requiring minimal additional firm capacity compared to the 44\% increase necessitated by the heuristic. Alongside each plan, reporting a cost limit (\UB) exceeded with probability no more than \(1-\alpha\) (i.e., Value-at-Risk bound via Cantelli's inequality) gives decision-makers a clear framework for evaluating the trade-off between cost and robustness.

Third, the robust plans mainly rely on changes in the type and location of generation and storage, with only moderate transmission expansion beyond the deterministic plan, to protect against line outages (Section~\ref{sec:cost_trade_off}). This finding contrasts with the predominant strategy in the reviewed robust and distributionally robust TEP literature, where transmission expansion is typically the only investment option to protect against contingencies \citep{5_Moreira2015_ARO, 9_Wu2018_robust, 18_Zhang2021_data, 8_Wang2022_dist, 4_Yuan2022_TEP}. The main reason for this is that in those studies, the generation portfolio is fixed exogenously rather than co-optimized with the network. When storage and network investment are co-optimized under contingencies at a smaller scale, similar substitution effects have been observed. For instance, \citet{huanca2024_TEP_storage} reported that incorporating storage into contingency-constrained transmission expansion reduces both transmission requirements and overall investment cost across three test systems. The results in this study extend this substitution effect to an inter-regional G\&TEP with long-duration hydrogen storage and duration-aware outages. Furthermore, our findings indicate that the location of storage and generation is as important as their quantity. The extent of this substitution, however, is likely sensitive to the assumed technology portfolio and network granularity (Section~\ref{sec:limitations}).

\subsection{Limitations and Future Work} \label{sec:limitations}

To focus on the effect of transmission outages, we adopted Assumptions~\ref{ass:single_year}--\ref{ass:corridor} in Section~\ref{sec:ProblemStatement}. With these assumptions, several limitations arise that act in both directions and make their net effect on the estimated outage risk unclear. In particular, the assumptions on uncertainty and outage processes may underestimate outage risk because they do not capture common-mode or cascading failures, nor the increased outage risk during extreme weather events. On the other hand, simplifications regarding corridor aggregation likely overstate the risk and make the robust plans more conservative. It may also understate the value of reinforcing a corridor with additional parallel circuits and thus favor generation and storage over transmission in the robust plans (Section~\ref{sec:cost_trade_off}). Furthermore, the 25-hour time step resolution cannot capture intraday VRE variability, so we excluded battery storage. We also neglected hydropower in our model. However, both batteries and hydropower could enhance local flexibility during outages and would likely affect how much firm generation and storage can substitute for transmission.

These simplifications suggest several avenues for future research. We recommend jointly treating weather, demand, and transmission uncertainty, and modeling generator outages as two-state Markov processes, similar to line outages. Failure probabilities can be linked to weather data to represent correlated and cascading failures \citep{32_panteli2017modeling, 31_allan2013reliability}, which the simulation function can accommodate without structural modifications. Additional extensions include circuit-level disaggregation of corridors and increased temporal resolution to incorporate short-duration battery storage alongside hydrogen storage. Applying the framework to larger systems will likely require scenario reduction, adaptive sampling, or more advanced memory and cut-aggregation strategies to remain tractable.

The findings above should therefore be treated as qualitative patterns from a stylized model, not as calibrated projections or investment recommendations for the modeled countries and policy decisions. This study's main contribution lies in its methodology and the general direction of the results, not in any specific cost or capacity estimates.

\section{Conclusion} \label{sec:conclusion}

This paper presented a two-stage robust optimization framework for inter-regional capacity expansion models that co-optimize generation, storage, and transmission investments. The framework treats transmission line outages as time-dependent scenarios that capture when and for how long lines fail (RQ\ref{rq:integration}). These scenarios are generated by Monte Carlo simulation, and the investment plan is hedged against the most severe of them. To solve the resulting model (RQ\ref{rq:solvability}), we proposed a Benders decomposition algorithm combined with a memory mechanism. In each iteration, the algorithm samples new scenarios, keeps the most critical ones, and solves the sub-problems in parallel. It uses Cantelli's inequality to bound the Value-at-Risk (VaR) of the system cost of each candidate plan at a chosen probability level $\alpha$. The selected plan is then evaluated on a large set of new, independently simulated scenarios, which gives more precise cost statistics. In a 25-node case study of Germany and its neighboring countries, the algorithm handled thousands of outage scenarios. After 300 iterations and an evaluation on 1000 simulated scenarios (14.7--17.9 hours of computation), the gap between its lower bound and the simulation-based VaR bound at $\alpha = 0.99$ was 0.28--0.94\%.

In the case study (RQ\ref{rq:comparison}), we compared the robust plans with two benchmarks, a deterministic plan that assumes all lines are always available and a reserve-margin heuristic. The heuristic is the deterministic model with the additional requirement that the firm capacity at each node is at least 1.1 times its peak demand. The deterministic plan had loss of load in up to 41\% of the simulated scenarios, with a mean of up to 79~\si{\GWh} and 1371~\si{\GWh} in the worst scenario. The heuristic nearly eliminated loss of load, but only by more than doubling gas capacity and increasing firm capacity by 44\%. The robust plans also nearly eliminated loss of load. Their mean adjusted system cost was 4.4--4.6\% below the heuristic and within 0.2\% of the deterministic plan, with the lowest cost volatility of the three approaches (Table~\ref{tab:results_cases}). The additional investment of 0.1--1.3\% over the deterministic plan was largely recovered through lower variable costs in the more severe outage cases. Rather than adding capacity, the robust approach changed the type and location of capacity. Part of the onshore wind capacity moved closer to demand, nuclear capacity increased by up to 13\%, hydrogen storage was spread over up to all 25 nodes, and transmission grew by only up to 6\% beyond the deterministic plan. 

These findings address the main gap identified in Section~\ref{sec:LiteratureReview}. The reviewed TEP and G\&TEP studies treat line outages as static snapshots and ignore when they occur and how long they last. Our results show that this simplification can hide important risks. Outage processes with the same expected unavailability led to similar mean loss of load but different worst cases, because rarer and longer outages caused larger events. Protecting against time-dependent outages also changed mainly the type and location of generation and storage, an option that most reviewed studies do not have because they fix the generation portfolio. Capacity expansion models for systems with high shares of wind and solar should therefore represent when and for how long lines fail, not only whether they fail.

\section*{Acknowledgments}

This work was funded by the Swedish Energy Agency under grant numbers P2023-01323 and P2022-00768. I.R. acknowledges support by the German Federal Ministry for Economic Affairs and Energy (BMWE) under Grant No.~03EI4083A (RESILIENT) jointly with the CETPartnership (\url{https://cetpartnership.eu/}) through the Joint Call 2022. As such, I.R. further acknowledges funding from the European Union's Horizon Europe research and innovation programme under grant agreement no.~101069750.

\section*{CRediT authorship contribution statement}

\textbf{Kamran Forghani:} Conceptualization, Data curation, Formal analysis, Investigation, Methodology, Software, Validation, Visualization, Writing -- original draft, Writing -- review \& editing.
\textbf{Lina Reichenberg:} Conceptualization, Funding acquisition, Project administration, Writing -- original draft, Writing -- review \& editing.
\textbf{Fredrik Hedenus:} Conceptualization, Funding acquisition, Project administration, Writing -- review \& editing.
\textbf{Iegor Riepin:} Conceptualization, Formal analysis, Investigation, Methodology, Validation, Writing -- review \& editing.

\section*{Declaration of competing interest}

The authors declare that they have no known competing financial interests or personal relationships that could have appeared to influence the work reported in this paper.

\section*{Data availability}

Data will be made available on request.

\bibliographystyle{elsarticle-num-names} 
\bibliography{cas-refs}

\onecolumn
\appendix

\counterwithin*{table}{section}
\counterwithin*{figure}{section}

\newpage
\section{Nomenclature} \label{app:nomenclature}

{\footnotesize
\begin{longtable}{@{}p{0.15\linewidth} p{0.65\linewidth} p{0.1\linewidth}@{}}
\caption{Model Nomenclature} \label{tab:nomenclature} \\
    \toprule
    \textbf{Symbol} & \textbf{Description} & \textbf{Unit} \\
    \midrule
    \endfirsthead
    \multicolumn{3}{c}%
    {{\bfseries \tablename\ \thetable{} -- continued from previous page}} \\
    \toprule
    \textbf{Symbol} & \textbf{Description} & \textbf{Unit} \\
    \midrule
    \endhead
    \midrule \multicolumn{3}{r}{{Continued on next page}} \\
    \endfoot
    \bottomrule
    \endlastfoot

    \multicolumn{3}{l}{\textbf{Sets and Indices}} \\
    $\setN, n$ & Set of, and index for, network nodes & - \\
    $n_{\text{0}}$ & Reference node of the voltage angles & - \\
    $\setL, l$ & Set of, and index for, transmission lines (inter-nodal corridors, Assumption~\ref{ass:corridor}), each defined by its sending and receiving nodes, $l = (m, n)$ & - \\
    $\setLac$ & Set of AC transmission lines ($\setLac \subseteq \setL$); the remaining lines are DC links & - \\
    $\setLin{n}$ & Set of lines with their receiving end connecting to node $n$ & - \\
    $\setLout{n}$ & Set of lines with their sending end connecting to node $n$ & - \\
    $\setT, t$ & Set of, and index for, time steps in the optimization horizon & - \\
    $\setTzero$ & Set of time steps including the initial state, $\setT \cup \{0\}$ & - \\
    $\setGgen, g$ & Set of, and index for, generation technologies & - \\
    $\setO$ & Set of all line status realizations, $\setO = \{0,1\}^{|\setL| \times |\setT|}$, with a generic element denoted by $\LINESTATUS$ & - \\
    $\setX$ & Feasible region for first-stage investment decisions, $\setX = \{\INVESTMENT \ge \bm{0} \mid \Ax\INVESTMENT \le \bx\}$ & - \\
    $\SampledScenarios$ & Set of outage realizations sampled independently from $\mathbb{P}$ in the scenario-based robust problem~\eqref{problem:sampled_robust} & - \\
    \\
    \multicolumn{3}{l}{\textbf{Parameters}} \\
    $\TIMESTEPDURATION$ & Duration of each time step & \si{\hour} \\
    $\OPWEIGHT$ & Weight of each time step in the annual operational cost, $\HOURSPERYEAR / |\setT|$ & \si{\hour} \\
    $\HOURSPERYEAR$ & Total hours per year, $\HOURSPERYEAR = 8760$ & \si{\hour} \\
    $\load{n}{t}$ & Electricity demand at node $n$ in time step $t$ & \si{\MW} \\
    $\cf{n}{g}{t}$ & Capacity factor of technology $g$ at node $n$ in time step $t$ & - \\
    $\maxcap{n}{g}, \mincap{n}{g}$ & Maximum and minimum installable capacity of technology $g$ at node $n$ & \si{\MW} \\
    $\maxcaptrans{l}, \mincaptrans{l}$ & Maximum and minimum installable capacity of transmission line $l$ & \si{\MW} \\
    $\ANNUALCOST{g}$ & Annualized investment and fixed O\&M cost of generation technology $g$ & \si{\Currency/MW/year} \\
    $\ANNUALCOST{\text{el}}, \ANNUALCOST{\text{fc}}$ & Annualized investment and fixed O\&M cost of electrolyzers (per MW of electricity input) and fuel cells (per MW of electricity output) & \si{\Currency/MW/year} \\
    $\ANNUALCOST{\text{H}_2}$ & Annualized investment and fixed O\&M cost of hydrogen storage & \si{\Currency/\MWh/year} \\
    $\ANNUALCOSTLine{l}$ & Annualized investment cost of transmission line $l$ & \si{\Currency/MW/year} \\
    $\VAROMCOST{g}$ & Variable cost of technology $g$ (variable O\&M plus fuel cost) & \si{\Currency/\MWh} \\
    $\COSTSHED$ & Penalty cost for voluntary load shedding & \si{\Currency/\MWh} \\
    $\COSTLOSS$ & Penalty cost for involuntary loss of load & \si{\Currency/\MWh} \\
    $\EMISSIONFACTOR{g}$ & Emission factor of technology $g$ & \si{\TonCOtwo/\MWh} \\
    $\CAPACITYCREDIT{g}$ & Capacity credit of technology $g$ & - \\
    $\RESERVEMARGIN$ & Reserve margin, the required ratio of local firm capacity to peak demand at each node & - \\
    $\MAXDEMANDCURTAILMENTRATIO$ & Maximum share of demand that can be shed voluntarily & - \\
    $\CARBONCAPINTENSITY$ & System-wide carbon intensity cap on served energy & \si{\TonCOtwo/\MWh} \\
    $\EFFICIENCYCHARGE, \EFFICIENCYDISCHARGE$ & Electrolyzer efficiency (hydrogen output per unit of electricity input) and fuel-cell efficiency (electricity output per unit of hydrogen input) & - \\
    $\Bline{l}$ & Susceptance of transmission line $l$ & p.u. \\
    $\Lline{l}$ & Length of transmission line $l$ & \si{\km} \\
    $\linestatus{l}{t}$ & Binary status of line $l$ in time step $t$ under realization $\LINESTATUS$ (1 = in service, 0 = out of service) & - \\
    $\FailRateBase$ & Base failure rate of transmission lines per unit length & \si{\per\hour\per\km} \\
    $\MTTR$ & Mean time to repair a failed transmission line & \si{\hour} \\
    $\FailRate{l}$ & Failure rate of line $l$, $\FailRate{l} = \FailRateBase \cdot \Lline{l}$ & \si{\per\hour} \\
    $\RepairRate$ & Repair rate of a failed line, $\RepairRate = 1/\MTTR$ & \si{\per\hour} \\
    $\ProbFail{l}$ & Probability that line $l$, if in service, fails from one time step to the next & - \\
    $\ProbRep{l}$ & Probability that line $l$, if out of service, is repaired from one time step to the next & - \\
    $\ProbInit{l}$ & Steady-state probability that line $l$ is in service & - \\
    $\ConfidenceLevel$ & Probability level of the Value-at-Risk, i.e., the required probability that the cost does not exceed the threshold (used in Cantelli's upper bound) & - \\
    $\noIterations$ & Number of Benders decomposition iterations in Phase~I & - \\
    $\noRandomScenarios$ & Number of new Monte Carlo outage scenarios generated per iteration & - \\
    $\noMemoryScenarios$ & Number of scenarios kept in the memory bank & - \\
    $\noSimScenarios$ & Number of out-of-sample outage scenarios used to evaluate a capacity plan (Phase~II) & - \\
    \\
    \multicolumn{3}{l}{\textbf{First-Stage (Investment) Variables}} \\
    $\Captech{n}{g}$ & Installed capacity of technology $g$ at node $n$ & \si{\MW} \\
    $\Captrans{l}$ & Installed capacity of transmission line $l$ & \si{\MW} \\
    $\CapEl{n}$ & Installed electrolyzer capacity (electricity input) at node $n$ & \si{\MW} \\
    $\CapFC{n}$ & Installed fuel-cell capacity (electricity output) at node $n$ & \si{\MW} \\
    $\CapHS{n}$ & Installed hydrogen storage capacity at node $n$ & \si{\MWh} \\
    $\TOC$ & Operational cost threshold, equal to the $\VaR_{\ConfidenceLevel}$ of the operational cost in~\eqref{problem:var} and to the worst-case operational cost over the sampled realizations in~\eqref{problem:sampled_robust} and in the master problem & \si{\Currency/year} \\
    $\INVESTMENT$ & Vector of first-stage investment decisions, $(\Captech{n}{g}, \CapEl{n}, \CapFC{n}, \CapHS{n}, \Captrans{l})$ & - \\
    \\
    \multicolumn{3}{l}{\textbf{Second-Stage (Operational) Variables}} \\
    $\Pgen{n}{g}{t}$ & Power dispatch of technology $g$ at node $n$ in time step $t$ & \si{\MW} \\
    $\Cin{n}{t}, \Cout{n}{t}$ & Electricity consumed by the electrolyzer and produced by the fuel cell at node $n$ in time step $t$ & \si{\MW} \\
    $\Estor{n}{t}$ & Hydrogen in storage at node $n$ at the end of time step $t$ & \si{\MWh} \\
    $\LineFlowFwd{l}{t}, \LineFlowBwd{l}{t}$ & Forward and backward power flow on line $l$ in time step $t$ & \si{\MW} \\
    $\ThetaVar{n}{t}$ & Voltage angle at node $n$ in time step $t$ (free variable) & rad \\
    $\Pshed{n}{t}$ & Voluntary load shedding at node $n$ in time step $t$ & \si{\MW} \\
    $\Ploss{n}{t}$ & Involuntary loss of load at node $n$ in time step $t$ & \si{\MW} \\
    $\OPERATION$ & Vector of second-stage operational decisions & - \\
    \\
    \multicolumn{3}{l}{\textbf{Algorithm Notation}} \\
    $\COSTX, \COSTy$ & Vectors of first-stage and second-stage cost coefficients & - \\
    $\Ax, \bx$ & Matrix and vector defining the first-stage feasible region, $\Ax\INVESTMENT \le \bx$ & - \\
    $\Wy(\LINESTATUS), \hy$ & Matrix of second-stage constraint coefficients, which depends on the line status through the power-flow equations, and vector of right-hand sides & - \\
    $\Ty{\LINESTATUS}$ & Matrix mapping a line status realization $\LINESTATUS$ to second-stage capacity limits & - \\
    $\Pi(\LINESTATUS)$ & Feasible set of the dual sub-problem for outage scenario $\LINESTATUS$ & - \\
    $\DualVariables$ & Vector of dual variables of the operational sub-problem & - \\
    $\LINESTATUS^{\text{base}}$ & Nominal scenario, in which all lines are always in service, $\bm{1}^{|\setL| \times |\setT|}$ & - \\
    $\Qx{\INVESTMENT}{\LINESTATUS}$ & Second-stage (recourse) operational cost of plan $\INVESTMENT$ under realization $\LINESTATUS$ & \si{\Currency/year} \\
    $\TC(\LINESTATUS)$ & Total system cost (investment and recourse cost) under realization $\LINESTATUS$ & \si{\Currency/year} \\
    $k$ & Iteration index of Phase~I & - \\
    $\overline{\INVESTMENT}$ & Candidate plan returned by the master problem in an iteration & - \\
    $\INVESTMENT^*$ & Incumbent plan, the candidate with the lowest in-loop upper bound & - \\
    $\RandomScenarios$ & Set of randomly generated outage scenarios in an iteration of Phase~I ($\noRandomScenarios$ scenarios) or in Phase~II ($\noSimScenarios$ scenarios) & - \\
    $\MemoryScenarios$ & Memory bank of the $\noMemoryScenarios$ most costly outage scenarios & - \\
    $\overline{\setO}$ & Set of outage scenarios evaluated in an iteration, $\overline{\setO} = \RandomScenarios \cup \MemoryScenarios$ & - \\
    $\LB$ & Lower bound on the optimal value of the scenario-based robust problem~\eqref{problem:sampled_robust} over the scenarios sampled so far & \si{\Currency/year} \\
    $\overline{\UB}$ & In-loop Cantelli bound of the candidate plan $\overline{\INVESTMENT}$, computed from the sample mean and variance of $\TC(\LINESTATUS)$ over $\overline{\setO}$ & \si{\Currency/year} \\
    $\UB$ & Upper bound, equal to the lowest $\overline{\UB}$ found so far in Phase~I and to the Cantelli bound on the $\VaR_{\ConfidenceLevel}$ of the total system cost of $\INVESTMENT^*$ estimated from the out-of-sample scenarios in Phase~II (reported) & \si{\Currency/year} \\
    $\Gap$ & Relative gap, $(\UB-\LB)/\UB \times 100$ & \si{\percent} \\
    $s$ & Index of a simulated scenario in Function~\ref{func:Simulate_Failures} & - \\
    $u$ & Random number drawn from $\RandU{0}{1}$ in Function~\ref{func:Simulate_Failures} & - \\
    \\
    \multicolumn{3}{l}{\textbf{Dual Variables (Operational Sub-problem)}} \\
    $\dualGen{n}{g}{t}$ & Dual variable of the generation dispatch limit~\eqref{eq:op_gen_limit} & - \\
    $\dualEl{n}{t}$ & Dual variable of the electrolyzer capacity limit~\eqref{eq:op_storage_power} & - \\
    $\dualFC{n}{t}$ & Dual variable of the fuel-cell capacity limit~\eqref{eq:op_storage_fc} & - \\
    $\dualHS{n}{t}$ & Dual variable of the hydrogen storage capacity limit~\eqref{eq:op_storage_energy} & - \\
    $\dualLine{l}{t}$ & Dual variable of the line flow capacity limit~\eqref{eq:op_line_cap} & - \\
    $\dualNode{n}{t}$ & Dual variable of the nodal energy balance~\eqref{eq:op_energy_balance} & - \\
    $\dualShed{n}{t}$ & Dual variable of the load shedding limit~\eqref{eq:op_shed_limit} & - \\
    $\dualUns{n}{t}$ & Dual variable of the unserved energy limit~\eqref{eq:op_unserved_limit} & - \\
    $\dualCO$ & Dual variable of the system carbon cap~\eqref{eq:op_carbon_cap} & - \\
    \\
    \multicolumn{3}{l}{\textbf{Results Metrics}} \\
    PoR & Price of robustness, the increase in annualized investment and fixed O\&M cost relative to the deterministic plan & \si{\Currency/year} \\
    $\VoRCost$ & Value of robustness, the reduction in mean variable cost (adjusted system cost excluding annualized investment and fixed O\&M costs) relative to the deterministic plan & \si{\Currency/year} \\
    $\VoRLoL$ & Value of robustness, the reduction in the highest loss of load across the simulated scenarios relative to the deterministic plan & \si{\GWh} \\
    \\
    \multicolumn{3}{l}{\textbf{Operators and Functions}} \\
    $\mathbb{P}$ & Probability measure on $\setO$ induced by the outage process & - \\
    $\VaR_{\ConfidenceLevel}(\cdot)$ & Value-at-Risk at probability level $\ConfidenceLevel$ ($\ConfidenceLevel$-quantile) & \si{\Currency/year} \\
    $\text{Mean}(\cdot), \text{Var}(\cdot)$ & Sample mean and sample variance over a set of scenarios & - \\
    $\arg\operatorname{top}_{\noMemoryScenarios}$ & Operator returning the $\noMemoryScenarios$ elements with the largest values & - \\
    $\Indicator{\cdot}$ & Indicator function, evaluating to $1$ if the condition holds and $0$ otherwise & - \\
    $\RandU{a}{b}$ & Random number drawn from the continuous uniform distribution on $[a, b]$ & - \\
\end{longtable}}

\section{Mathematical Models}

\subsection{Master Problem} \label{app:master_model}

The master problem (MP) in our decomposition algorithm is a deterministic model formulated as follows.  

{\small
\begin{alignat}{3}
    \textbf{MP:} \quad \min \quad & \sum_{n \in \setN} \sum_{g \in \setGgen} \ANNUALCOST{g} \cdot \Captech{n}{g} + \sum_{n \in \setN} \left( \ANNUALCOST{\text{el}} \cdot \CapEl{n} + \ANNUALCOST{\text{fc}} \cdot \CapFC{n} + \ANNUALCOST{\text{H}_2} \cdot \CapHS{n} \right) + \sum_{l \in \setL} \ANNUALCOSTLine{l} \cdot \Captrans{l} + \TOC \label{eq:system_cost}\\
    \text{s.t.:} \quad
    & \mincap{n}{g} \le \Captech{n}{g} \leq \maxcap{n}{g} && \forall n \in \setN, \; g \in \setGgen \label{eq:inv_res_potential}\\
    & \mincaptrans{l} \leq \Captrans{l} \leq \maxcaptrans{l} && \forall l \in \setL \label{eq:inv_trans_potential}\\
    & \frac{\TIMESTEPDURATION}{\EFFICIENCYDISCHARGE} \cdot \CapFC{n} \le \CapHS{n} && \forall n \in \setN \label{eq:inv_h2_fc_energy}\\
    & \EFFICIENCYCHARGE \cdot \TIMESTEPDURATION \cdot \CapEl{n} \le \CapHS{n} && \forall n \in \setN \label{eq:inv_h2_el_energy}\\
    & \CapFC{n} \le \EFFICIENCYCHARGE \cdot \EFFICIENCYDISCHARGE \cdot |\setT| \cdot \CapEl{n} && \forall n \in \setN \label{eq:inv_h2_fc_el}\\
    & \TOC \ge \Qx{\INVESTMENT}{\LINESTATUS^\text{base}} && \text{where } \LINESTATUS^\text{base} = \bm{1}^{|\setL| \times |\setT|} \label{eq:initial_solution}\\
    & \text{All first-stage variables} \geq 0 \nonumber
\end{alignat}}

In the model above, the objective function~\eqref{eq:system_cost} minimizes the system cost, comprising the annualized investment and fixed O\&M costs of generation, hydrogen storage, and transmission, and the variable system cost, $\TOC$. Through the optimality cuts~\eqref{eq:benders_cut} added in Algorithm~\ref{alg:Benders_Decomposition}, $\TOC$ approximates from below the worst-case operational cost over the sampled outage scenarios, as in Problem~\eqref{problem:sampled_robust}. Constraint sets~\eqref{eq:inv_res_potential} and \eqref{eq:inv_trans_potential}, respectively, ensure that the installed capacity of the technologies and transmission lines remains within predefined limits.

Constraints~\eqref{eq:inv_h2_fc_energy}--\eqref{eq:inv_h2_fc_el} link the sizes of the electrolyzer, fuel cell, and hydrogen storage at each node. Without them, the MP could build one component without the others, for example, a fuel cell without a storage unit or electrolyzer to supply it, because the MP does not see how the storage is operated. The constraints only rule out capacity that could never be used in any scenario, so they do not change the optimal solution. Constraint~\eqref{eq:inv_h2_fc_energy} requires that a full hydrogen storage unit can run the fuel cell at full power for one time step. Constraint~\eqref{eq:inv_h2_el_energy} requires that an empty storage unit can take in the full output of the electrolyzer for one time step. Constraint~\eqref{eq:inv_h2_fc_el} requires that the fuel cell does not deliver more electricity in one time step than the electrolyzer can supply over the whole year, after both conversion losses. This last constraint holds because hydrogen enters the store only through the electrolyzer and the storage level is the same at the start and the end of the year (Constraints~\eqref{eq:op_storage_balance} and~\eqref{eq:op_storage_cyclic}).

Constraint~\eqref{eq:initial_solution} adds the operation of the nominal scenario, in which all lines are always available, $\LINESTATUS^{\text{base}} = \bm{1}^{|\setL| \times |\setT|}$, to the MP, with all its operational variables and constraints (\ref{app:subproblem_model}). The operational cost $\TOC$ of the MP is therefore at least the cost of this scenario, $\Qx{\INVESTMENT}{\LINESTATUS^{\text{base}}}$. Without this constraint, the MP would see no operational cost in the first iteration and build no generation capacity at all, which slows down convergence. With it, the first MP is the deterministic model, so the algorithm starts from the deterministic plan and a good lower bound.
\subsection{Second-Stage Sub-Problem} \label{app:subproblem_model}

For a fixed investment plan $\INVESTMENT$ (i.e., $\Captech{n}{g}, \CapEl{n}, \CapFC{n}, \CapHS{n}, \Captrans{l}$) and a specific transmission line outage scenario $\LINESTATUS$ (i.e., $\linestatus{l}{t}$), the second-stage sub-problem is given below. The dual variable corresponding to each constraint is shown in brackets in front of it. Note that dual variables associated with constraints having a zero right-hand side are omitted, as they do not contribute to cut generation (see \eqref{eq:benders_cut}).\\

\begin{adjustbox}{width=\textwidth}
$\begin{alignedat}{3}
    \textbf{SP:} \quad \Qx{\INVESTMENT}{\LINESTATUS} = \min \quad & \sum_{n \in \setN} \sum_{g \in \setGgen} \sum_{t \in \setT} \OPWEIGHT \cdot \VAROMCOST{g} \cdot \Pgen{n}{g}{t} + \sum_{n \in \setN} \sum_{t \in \setT} \OPWEIGHT \cdot \left( \COSTSHED \cdot \Pshed{n}{t} + \COSTLOSS \cdot \Ploss{n}{t} \right) \qquad \span \span \label{eq:toc_def} \\
    \text{s.t.:} \quad
    & \Pgen{n}{g}{t} \leq \cf{n}{g}{t} \cdot \Captech{n}{g} && \forall n \in \setN, \; g \in \setGgen, \; t \in \setT && \qquad [\dualGen{n}{g}{t}] \label{eq:op_gen_limit}\\
    & \Cin{n}{t} \leq \CapEl{n} && \forall n \in \setN, \; t \in \setT && \qquad [\dualEl{n}{t}] \label{eq:op_storage_power} \\
    & \Cout{n}{t} \leq \CapFC{n} && \forall n \in \setN, \; t \in \setT && \qquad [\dualFC{n}{t}] \label{eq:op_storage_fc} \\
    & \Estor{n}{t} \leq \CapHS{n} && \forall n \in \setN, \; t \in \setTzero && \qquad [\dualHS{n}{t}] \label{eq:op_storage_energy}\\
    & \LineFlowFwd{l}{t} + \LineFlowBwd{l}{t} \leq \linestatus{l}{t} \cdot \Captrans{l} && \forall l \in \setL, \; t \in \setT && \qquad [\dualLine{l}{t}] \label{eq:op_line_cap}\\
    & \LineFlowFwd{l}{t} - \LineFlowBwd{l}{t} - \Bline{l} \cdot \linestatus{l}{t} (\ThetaVar{m}{t} - \ThetaVar{n}{t}) = 0 && \forall (m, n) = l \in \setLac, \; t \in \setT \label{eq:op_dc_flow}\\
    & \ThetaVar{n_{\text{0}}}{t} = 0 && \forall t \in \setT \label{eq:op_ref_angle}\\
    & \Estor{n}{t-1}  + \TIMESTEPDURATION \left(\Cin{n}{t} \cdot \EFFICIENCYCHARGE - \frac{\Cout{n}{t}}{\EFFICIENCYDISCHARGE}\right) - \Estor{n}{t} = 0 && \forall n \in \setN, \; t \in \setT \label{eq:op_storage_balance} \\
    & \Estor{n}{|\setT|} - \Estor{n}{0} = 0 && \forall n \in \setN \label{eq:op_storage_cyclic}\\
    & \sum_{g \in \setGgen} \Pgen{n}{g}{t} + \Cout{n}{t} - \Cin{n}{t} + \sum_{l \in \setLin{n}} (\LineFlowFwd{l}{t} - \LineFlowBwd{l}{t}) \nonumber \\ 
    & \qquad \qquad + \sum_{l \in \setLout{n}} (\LineFlowBwd{l}{t} - \LineFlowFwd{l}{t}) + \Pshed{n}{t} + \Ploss{n}{t} = \load{n}{t}   && \forall n \in \setN, \; t \in \setT && \qquad [\dualNode{n}{t}] \label{eq:op_energy_balance} \\
    & \Pshed{n}{t} \leq \MAXDEMANDCURTAILMENTRATIO \cdot \load{n}{t} && \forall n \in \setN, \; t \in \setT && \qquad [\dualShed{n}{t}] \label{eq:op_shed_limit}\\
    & \Pshed{n}{t} + \Ploss{n}{t} \leq \load{n}{t} && \forall n \in \setN, \; t \in \setT && \qquad [\dualUns{n}{t}] \label{eq:op_unserved_limit}\\
    & \sum_{n \in \setN} \sum_{g \in \setGgen} \sum_{t \in \setT} \EMISSIONFACTOR{g} \cdot \Pgen{n}{g}{t} \leq \CARBONCAPINTENSITY \sum_{n \in \setN} \sum_{t \in \setT} \left( \load{n}{t} - \Pshed{n}{t} - \Ploss{n}{t} \right) && && \qquad [\dualCO] \label{eq:op_carbon_cap}\\
    & \text{All second-stage variables, except $\ThetaVar{n}{t}$, are nonnegative} \nonumber 
\end{alignedat}$
\end{adjustbox}\\

In the model above, for a given investment plan $\INVESTMENT$ and outage scenario $\LINESTATUS$, objective function~\eqref{eq:toc_def} minimizes variable cost, consisting of variable generation costs and penalties for voluntary load shedding and involuntary loss of load. Each time step is weighted by $\OPWEIGHT = \HOURSPERYEAR / |\setT|$, which scales the operation over the modeled time steps to a full year. Constraint set~\eqref{eq:op_gen_limit} states that the generated electricity in each node at each time step using a generation technology is restricted by the corresponding capacity factor. For variable renewable sources such as wind and solar, $\cf{n}{g}{t}$ represents the capacity factor time series, taking values in $[0, 1]$, whereas for dispatchable sources, such as nuclear or gas, it is set to 1. Constraint sets \eqref{eq:op_storage_power} and \eqref{eq:op_storage_fc} limit the electricity consumed by the electrolyzer and the electricity produced by the fuel cell to their installed capacities, respectively. Constraint set \eqref{eq:op_storage_energy} ensures that the hydrogen stored at each node, measured in energy units, does not exceed the installed storage capacity. Constraint set \eqref{eq:op_line_cap} states that the sum of forward and backward power flow through the transmission lines cannot exceed the installed transmission capacity. Note that the constraint is written in terms of the combined directional flows, since at optimality, for any time step, at least one of these flows will be zero.

Constraint set \eqref{eq:op_dc_flow} is a DC approximation of Kirchhoff’s Voltage Law (KVL) for power flow through AC transmission lines\footnote{In this linearization, voltage magnitudes are assumed to be close to unity, reactive power flows and conductance effects are neglected, and the real power flow on each line is proportional to the difference in voltage angles at its terminals, which results in a set of linear constraints relating power flows to phase angle differences around network cycles, see~\citet{38_Neumann2022_assessments} for more details.}. Constraint set \eqref{eq:op_ref_angle} removes the arbitrary degree of freedom in voltage angles and ensures a unique solution for the DC power flow equations. Constraint set~\eqref{eq:op_storage_balance} is the storage balance constraint that determines the hydrogen in storage at the end of each time step from the storage level in the previous time step, the hydrogen produced by the electrolyzer, and the hydrogen consumed by the fuel cell. Constraint set~\eqref{eq:op_storage_cyclic} ensures that the storage balance constraint remains closed-loop by connecting the first and last time steps. Constraint set~\eqref{eq:op_energy_balance} ensures nodal supply-demand balance by enforcing Kirchhoff’s Current Law (KCL). It requires that the net sum of generation, fuel-cell output minus electrolyzer consumption, line flows, load shedding, and loss of load matches the local electricity demand at each node for every time step. Constraint set~\eqref{eq:op_shed_limit} ensures that total voluntary load shedding at each node and time step does not exceed a certain fraction, $\MAXDEMANDCURTAILMENTRATIO$, of the electricity load in that period. Constraint set~\eqref{eq:op_unserved_limit} ensures that load shedding and loss of load together do not exceed the local demand. Finally, Constraint~\eqref{eq:op_carbon_cap} limits total emissions to $\CARBONCAPINTENSITY$ times the satisfied load.

\subsection{Optimality Cut} \label{app:optimality_cut}

In the decomposition algorithm, after solving the recourse sub-problems, for every $\LINESTATUS \in \overline{\setO}$, the following optimality cut is added to the master problem (see Line~\ref{line:add_cuts} in Algorithm~\ref{alg:Benders_Decomposition}). 

\begin{align}
    \TOC & \ge \sum_{n \in \setN} \sum_{g \in \setGgen} \sum_{t \in \setT} \dualGen{n}{g}{t} \cdot \cf{n}{g}{t} \cdot \Captech{n}{g}
    + \sum_{n \in \setN} \sum_{t \in \setT} \left( \dualEl{n}{t} \cdot \CapEl{n} + \dualFC{n}{t} \cdot \CapFC{n} \right)
    + \sum_{n \in \setN} \sum_{t \in \setTzero} \dualHS{n}{t} \cdot \CapHS{n}
    + \sum_{l \in \setL} \sum_{t \in \setT} \dualLine{l}{t} \cdot \linestatus{l}{t} \cdot \Captrans{l} \nonumber \\
    & + \sum_{n \in \setN} \sum_{t \in \setT} \left( \dualNode{n}{t} + \MAXDEMANDCURTAILMENTRATIO \cdot \dualShed{n}{t} + \dualUns{n}{t} \right) \cdot \load{n}{t}
    + \dualCO \cdot \CARBONCAPINTENSITY \sum_{n \in \setN} \sum_{t \in \setT} \load{n}{t} \label{eq:benders_cut}
\end{align}

where the coefficients, denoted by $\pi$, are the optimal dual variables (shadow prices) associated with the corresponding primal constraints in \eqref{eq:op_gen_limit}--\eqref{eq:op_line_cap} and \eqref{eq:op_energy_balance}--\eqref{eq:op_carbon_cap} after solving the sub-problem for outage scenario $\LINESTATUS$. The cut contains one term for every sub-problem constraint with a nonzero right-hand side. By linear programming duality, the right-hand side of~\eqref{eq:benders_cut} equals $\Qx{\overline{\INVESTMENT}}{\LINESTATUS}$ at the candidate plan $\overline{\INVESTMENT}$ at which the sub-problem was solved; we verify this equality numerically for every generated cut as a check that the cut is complete.
\subsection{Reserve Margin Constraint} \label{app:reserve_margin}

As a heuristic approach, we add Constraint~\eqref{eq:inv_res_margin} to the master problem (MP) presented in \ref{app:master_model} to deal with the transmission grid uncertainty. This constraint ensures that the total firm capacity at each node meets a certain fraction, $\RESERVEMARGIN$, of the peak electricity demand through the planning horizon. Only generation capacity counts toward firm capacity; hydrogen fuel cells are not credited. Higher $\RESERVEMARGIN$ would reduce reliance on the transmission grid; thus, as a heuristic, it would increase robustness at the expense of higher system cost. 

\begin{equation}
    \sum_{g \in \setGgen} \CAPACITYCREDIT{g} \cdot \Captech{n}{g} \geq \RESERVEMARGIN \cdot \max_{t \in \setT} \left\{ \load{n}{t} \right\} \quad \forall n \in \setN \label{eq:inv_res_margin}
\end{equation}
\section{Simulation of Transmission Line Outages} \label{app:simulation}

Function~\ref{func:Simulate_Failures} generates the set of outage scenarios $\RandomScenarios$ of Section~\ref{sec:meth_simulation} by Monte Carlo simulation. Each scenario is a matrix $\LINESTATUS$ with one row for each line $l \in \setL$ and one column for each time step $t \in \setT$. Its entry $\linestatus{l}{t}$ is $1$ if line $l$ is in service in time step $t$ and $0$ if it is out of service.

The first part of the function (Lines~\ref{alg2:line:repair_rate}--\ref{alg2:line:precalc_end}) computes the rates and probabilities of each line $l$. They are the same in all scenarios, so they are computed only once. Line~\ref{alg2:line:repair_rate} computes the repair rate and Line~\ref{alg2:line:fail_rate} the failure rate of line $l$. Lines~\ref{alg2:line:fail_prob} and~\ref{alg2:line:rep_prob} compute the per-step failure and repair probabilities of Eq.~\eqref{eq:transition_probabilities}, and Line~\ref{alg2:line:steady_state} computes the steady-state probability that line $l$ is in service.

The second part (Lines~\ref{alg2:line:mc_loop_start}--\ref{alg2:line:mc_loop_end}) creates the scenarios one by one. Line~\ref{alg2:line:init_matrix} creates the scenario matrix for $s$. For each line $l$, Line~\ref{alg2:line:init_state} sets its status in the first time step. For each later time step, Line~\ref{alg2:line:draw_u} draws a random number $u$ between $0$ and $1$, and Line~\ref{alg2:line:transition} sets the new status of line $l$ based on its status in the previous time step. If line $l$ was in service ($\linestatus{l}{t-1} = 1$), only the first term of Line~\ref{alg2:line:transition} matters, and line $l$ stays in service if $u \ge \ProbFail{l}$ and fails otherwise. If line $l$ was out of service ($\linestatus{l}{t-1} = 0$), only the second term matters, and line $l$ is repaired if $u < \ProbRep{l}$ and stays out of service otherwise. Line~\ref{alg2:line:add_scenario} adds the finished matrix to the set $\RandomScenarios$, which is empty at the start (Line~\ref{alg2:line:init}), and Line~\ref{alg2:line:return} returns all scenarios.

\begin{function}[htbp]
	\caption{Simulate($\noRandomScenarios, \FailRateBase, \MTTR$)} \label{func:Simulate_Failures}
	\LinesNumbered
	\DontPrintSemicolon

	$\RandomScenarios \gets \varnothing$ \label{alg2:line:init} \tcp*{empty set of scenarios}

	\tcp{part 1: probabilities of each line $l$ (the same in all scenarios)}
	$\RepairRate \gets 1 / \MTTR$ \label{alg2:line:repair_rate} \tcp*{repair rate}
	\For{$l \in \setL$ \label{alg2:line:precalc_start}}{
		$\FailRate{l} \gets \FailRateBase \cdot \Lline{l}$ \label{alg2:line:fail_rate} \tcp*{failure rate, grows with line length}
		$\ProbFail{l} \gets \frac{\FailRate{l}}{\FailRate{l}+\RepairRate} \big(1 - \exp\big(-(\FailRate{l} + \RepairRate) \cdot \TIMESTEPDURATION\big)\big)$ \label{alg2:line:fail_prob} \tcp*{failure probability per time step}
		$\ProbRep{l} \gets \frac{\RepairRate}{\FailRate{l}+\RepairRate} \big(1 - \exp\big(-(\FailRate{l} + \RepairRate) \cdot \TIMESTEPDURATION\big)\big)$ \label{alg2:line:rep_prob} \tcp*{repair probability per time step}
		$\ProbInit{l} \gets \frac{\RepairRate}{\FailRate{l}+\RepairRate}$ \label{alg2:line:steady_state} \tcp*{steady-state probability of being in service}
	} \label{alg2:line:precalc_end}

	\tcp{part 2: create the scenarios}
	\For{$s \gets 1$ \KwTo $\noRandomScenarios$ \label{alg2:line:mc_loop_start}}{
		$\LINESTATUS \gets \bm{0}^{|\setL| \times |\setT|}$ \label{alg2:line:init_matrix} \tcp*{matrix of scenario $s$}
		\For{$l \in \setL$ \label{alg2:line:line_loop_start}}{
			$\linestatus{l}{1} \gets \Indicator{\RandU{0}{1} < \ProbInit{l}}$ \label{alg2:line:init_state} \tcp*{status at $t = 1$}
			\For{$t \gets 2$ \KwTo $|\setT|$ \label{alg2:line:temporal_loop_start}}{
				$u \gets \RandU{0}{1}$ \label{alg2:line:draw_u} \tcp*{random number between $0$ and $1$}
				$\linestatus{l}{t} \gets \linestatus{l}{t-1} \cdot \Indicator{u \ge \ProbFail{l}} + (1 - \linestatus{l}{t-1}) \cdot \Indicator{u < \ProbRep{l}}$ \label{alg2:line:transition} \tcp*{new status}
			} \label{alg2:line:temporal_loop_end}
		} \label{alg2:line:line_loop_end}
		$\RandomScenarios \gets \RandomScenarios \cup \{\LINESTATUS\}$ \label{alg2:line:add_scenario} \tcp*{add scenario $s$ to the set}
	} \label{alg2:line:mc_loop_end}
	\Return{$\RandomScenarios$} \label{alg2:line:return} \tcp*{all scenarios}
\end{function}


\section{Further Details of the Decomposition Algorithm} \label{app:benders_details}

\begin{figure}[htbp]
    \centering
    \resizebox{0.7\textwidth}{!}{%
        \begin{tikzpicture}[
    scale=1, transform shape,
    node distance=0.55cm and 1cm,
    startstop/.style={ellipse, minimum width=1.5cm, minimum height=0.9cm, align=center, draw=black, fill=red!10, line width=1pt, font=\bfseries},
    process/.style={rectangle, minimum width=2cm, minimum height=0.8cm, align=center, draw=black, fill=blue!5, line width=0.8pt},
    splitproc/.style={
        rectangle split,
        rectangle split parts=2,
        rectangle split part fill={blue!5, blue!10},
        minimum width=2cm,
        align=center,
        draw=black,
        line width=0.8pt
    },
    splitsubproc/.style={
        rectangle split,
        rectangle split parts=2,
        rectangle split part fill={orange!10, orange!15},
        dashed,
        minimum width=2cm,
        align=center,
        draw=black,
        line width=0.8pt
    },
    evalproc/.style={
        rectangle split,
        rectangle split parts=2,
        rectangle split part fill={green!6, green!12},
        minimum width=2cm,
        align=center,
        draw=black,
        line width=0.8pt
    },
    decision/.style={diamond, aspect=3, minimum width=3.5cm, minimum height=1cm, inner sep=1pt, align=center, draw=black, fill=yellow!15, line width=0.8pt},
    arrow/.style={thick, ->, >=Stealth},
    colA/.style={minimum width=5.65cm},
    colB/.style={minimum width=6.25cm},
    colC/.style={minimum width=4.85cm},
    phaseframe/.style={rounded corners=6pt, line width=0.8pt, inner xsep=10pt, inner ysep=10pt},
    phasetitle/.style={anchor=south east, inner sep=3pt, font=\bfseries}
    ]


    \node (start) [startstop] {Start};

    \node (init) [splitproc, colA, below=of start] {
        \textbf{Initialization} \\
        $k \gets 1$, $\LB \gets -\infty$, $\UB \gets \infty$ \\
        Fill memory $\MemoryScenarios$ with \\
        $\noMemoryScenarios$ simulated scenarios \\
        Initialize MP constraints $\setX$ \\
        with nominal grid status $\LINESTATUS^\text{base} = \bm{1}$
        \nodepart{two}
        \textbf{Output:} $(\setX, \MemoryScenarios, \LB, \UB)$
    };

    \node (check_iter) [decision, below=of init] {Is $k \le \noIterations$?};

    \node (master) [splitproc, colA, below=of check_iter] {
        \textbf{Solve Master Problem (MP)} \\
        Obtain lower bound $\LB$ \\
        and new investment plan $\overline{\INVESTMENT}$ \\
        by solving MP over $\setX$
        \nodepart{two}
        \textbf{Output:} $(\overline{\INVESTMENT}, \LB)$
    };

    \node (sim) [splitproc, colA, below=of master] {
        \textbf{Monte Carlo Simulation} \\
        Generate $\noRandomScenarios$ simulated scenarios $\RandomScenarios$ \\
        $\overline{\setO} \gets \RandomScenarios \cup \MemoryScenarios$
        \nodepart{two}
        \textbf{Output:} $(\RandomScenarios, \overline{\setO})$
    };

    \node (duals) [splitsubproc, colA, below=of sim] {
        \textbf{Solve Subproblems in Parallel} \\
        For each $\LINESTATUS \in \overline{\setO}$, obtain the system \\
        cost $\TC(\LINESTATUS)$ and dual values $\overline{\DualVariables}(\LINESTATUS)$
        \nodepart{two}
        \textbf{Output:} $(\overline{\DualVariables}(\LINESTATUS), \TC(\LINESTATUS))$
    };

    \node (memory) [splitproc, colA, below=of duals] {
        \textbf{Update Memory} \\
        Update $\MemoryScenarios$ considering $\TC(\LINESTATUS)$
        \nodepart{two}
        \textbf{Output:} $\MemoryScenarios$
    };


    \node (cuts) [splitproc, colB, anchor=north west] at ([xshift=1.0cm]master.north east) {
        \textbf{Cut Generation} \\
        Add Benders cuts to MP $\forall \LINESTATUS \in \overline{\setO}$: \\
        $\setX \gets \setX \cap \{ \TOC \ge \overline{\DualVariables}(\LINESTATUS)^T (\hy - \Ty{\LINESTATUS} \INVESTMENT) \}$
        \nodepart{two}
        \textbf{Output:} $\setX$
    };

    \node (calc_ub) [splitproc, colB, below=of cuts] {
        \textbf{In-Loop Cantelli Upper Bound} \\
        Calculate current upper bound $\overline{\UB}$ \\
        from $\TC(\LINESTATUS)$, $\LINESTATUS \in \overline{\setO}$
        \nodepart{two}
        \textbf{Output:} $\overline{\UB}$
    };

    \node (check_best) [decision, below=of calc_ub] {Is $\overline{\UB} < \UB$?};

    \node (update_best) [splitproc, right=0.8cm of check_best] {
        \textbf{Update Incumbent} \\
        $\UB \gets \overline{\UB}$, $\INVESTMENT^* \gets \overline{\INVESTMENT}$
        \nodepart{two}
        \textbf{Output:} $(\UB , \INVESTMENT^*)$
    };

    \node (calc_gap) [splitproc, colB, below=of check_best] {
        \textbf{Calculate In-Loop Gap} \\
        $\Gap \gets \frac{\UB - \LB}{\UB} \times 100\%$
        \nodepart{two}
        \textbf{Output:} $\Gap$
    };

    \node (increment) [process, below=of calc_gap] {$k \gets k + 1$};

    \coordinate (laneR) at ([xshift=0.5cm]update_best.east);
    \coordinate (laneL) at ([xshift=-0.65cm]init.west);
    \coordinate (corridor) at ($(memory.east)!0.5!(memory.east -| cuts.west)$);

    \begin{scope}[on background layer]
        \node (phase1) [phaseframe, draw=blue!50!black, fill=blue!2,
                        fit=(init) (memory) (cuts) (update_best) (increment) (laneR)] {};
    \end{scope}
    \node [phasetitle, text=blue!50!black] at (phase1.north east)
        {Phase I: Memory-augmented Benders decomposition};


    \node (eval_sim) [evalproc, colC, anchor=north west] at ([xshift=10pt, yshift=-1.25cm]phase1.south west) {
        \textbf{Out-of-Sample Simulation} \\
        Generate $\noSimScenarios$ new scenarios \\
        $\RandomScenarios$, independent of all \\
        scenarios used in Phase I
        \nodepart{two}
        \textbf{Output:} $\RandomScenarios$
    };

    \node (eval_sp) [splitsubproc, colC, right=0.55cm of eval_sim] {
        \textbf{Evaluate $\INVESTMENT^*$ in Parallel} \\
        For each $\LINESTATUS \in \RandomScenarios$, solve the \\
        subproblem for fixed $\INVESTMENT^*$ \\
        $\TC(\LINESTATUS) = \COSTX^T \INVESTMENT^* + \Qx{\INVESTMENT^*}{\LINESTATUS}$
        \nodepart{two}
        \textbf{Output:} $\TC(\LINESTATUS)$
    };

    \node (eval_ub) [evalproc, colC, right=0.55cm of eval_sp] {
        \textbf{Final Bound \& Gap} \\
        $\UB \gets$ Cantelli bound from \\
        sample moments of $\TC(\LINESTATUS)$ \\
        $\Gap \gets \frac{\UB - \LB}{\UB} \times 100\%$
        \nodepart{two}
        \textbf{Output:} $(\LB, \UB, \Gap, \INVESTMENT^*)$
    };

    \begin{scope}[on background layer]
        \node (phase2) [phaseframe, draw=green!40!black, fill=green!2,
                        fit=(eval_sim) (eval_sp) (eval_ub)] {};
    \end{scope}
    \node [phasetitle, text=green!40!black] at (phase2.north east)
        {Phase II: Out-of-sample evaluation of the incumbent plan $\INVESTMENT^*$};

    \node (stop) [startstop, right=0.5cm of phase2.east] {End};


    \draw [arrow] (start) -- (init);
    \draw [arrow] (init) -- (check_iter);
    \draw [arrow] (check_iter) -- node[right] {Yes} (master);
    \draw [arrow] (master) -- (sim);
    \draw [arrow] (sim) -- (duals);
    \draw [arrow] (duals) -- (memory);

    \draw [arrow] (memory.east) -- (memory.east -| corridor) |- (cuts.west);

    \draw [arrow] (cuts) -- (calc_ub);
    \draw [arrow] (calc_ub) -- (check_best);
    \draw [arrow] (check_best) -- node[above] {Yes} (update_best);
    \draw [arrow] (check_best) -- node[right] {No} (calc_gap);
    \draw [arrow] (update_best.south) |- (calc_gap.east);
    \draw [arrow] (calc_gap) -- (increment);

    \draw [arrow] (increment.east) -- (increment.east -| laneR) |- (check_iter.east);

    \draw [arrow] (check_iter.west) -- node[above, pos=0.3] {No} (check_iter.west -| laneL) |- (eval_sim.west);
    \draw [arrow] (eval_sim) -- (eval_sp);
    \draw [arrow] (eval_sp) -- (eval_ub);
    \draw [arrow] (eval_ub.east) -- (eval_ub.east -| stop.west);

\end{tikzpicture}
    }
    \caption{Flowchart of the proposed algorithm. Phase~I chooses the incumbent plan with the adaptive, memory-augmented Benders decomposition. Phase~II evaluates this plan on new out-of-sample scenarios. Dashed boxes are solved in parallel.}
    \label{fig:benders_flowchart}
\end{figure}
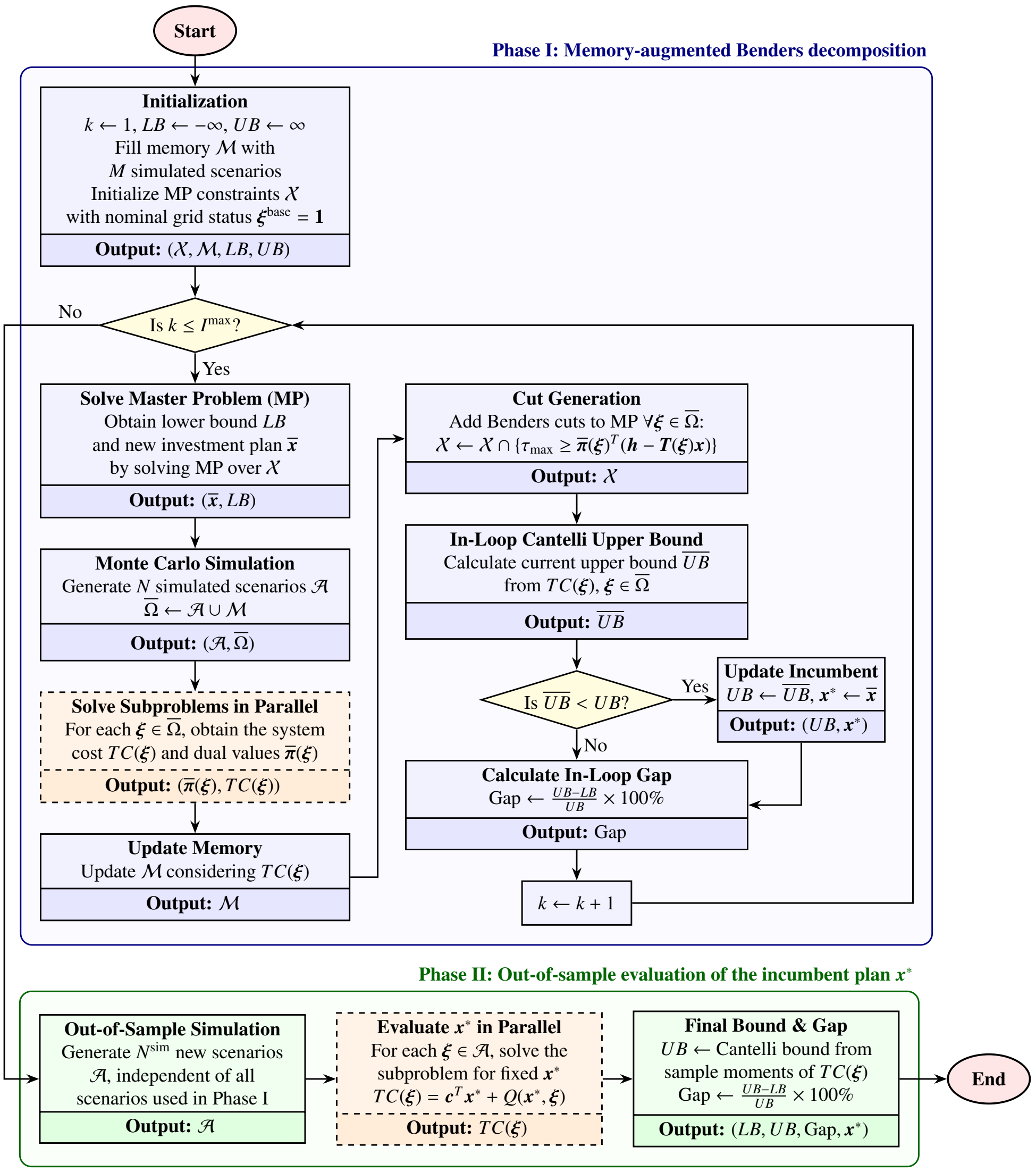

Figure~\ref{fig:benders_flowchart} shows Algorithm~\ref{alg:Benders_Decomposition} as a flowchart. The rest of this appendix explains the properties of the in-loop upper bound used in Phase~I (Section~\ref{sec:benders_decomposition}).

In each iteration, we compute the upper bound $\overline{\UB}$ by inserting the sample mean and sample standard deviation of the total system cost over the evaluated scenarios $\overline{\setO}$ into the Cantelli bound~\eqref{eq:cantelli} (Line~\ref{line:cantelli} of Algorithm~\ref{alg:Benders_Decomposition}). We use this sample-based estimate for two purposes: to select the reported plan and to check convergence. Therefore, three properties of this estimate matter when we interpret the results.

First, $\overline{\UB}$ uses estimated moments, not the true moments. So \eqref{eq:cantelli} holds only approximately for $\overline{\UB}$, and we do not attach any formal statistical confidence to $\UB$.

Second, $\overline{\setO}$ combines two kinds of scenarios. The random batch $\RandomScenarios$ is drawn after $\overline{\INVESTMENT}$ is fixed. Therefore, it gives an unbiased sample of the cost of $\overline{\INVESTMENT}$, as in standard Monte Carlo upper-bounding procedures \citep{mak1999monte}. In contrast, the memory scenarios $\MemoryScenarios$ are kept because they were the most costly realizations observed so far (Section~\ref{sec:benders_decomposition}). Including them over-represents the upper tail of the cost distribution. This usually raises the sample mean and, as long as the memory scenarios remain costly under $\overline{\INVESTMENT}$, also raises the sample standard deviation. We accept this bias on purpose, because it makes $\overline{\UB}$ larger and therefore more conservative.

Third, the incumbent plan is the candidate with the lowest $\overline{\UB}$ among many noisy estimates. Therefore, its in-loop $\UB$ may be slightly optimistic. This partly offsets the conservative effect of the memory scenarios.

For these reasons, we use the in-loop bound only to select the incumbent plan and to check convergence. The reported $\UB$ is computed in Phase~II (Section~\ref{sec:benders_decomposition}). This avoids the second and third effects and reduces the sampling error of the first.

\section{Input Data and System Parameters} \label{app:data}

The network topology, transmission capacities and costs, hourly demand, and renewable capacity factors and potentials come from a 25-node network generated with \texttt{PyPSA-Eur} \citep{horsch2018pypsa} (v2025.07.0) and based on the OpenStreetMap transmission grid \citep{xiong2025modelling}. These data refer to the weather year 2013, with the historical demand of that year and capacity factors computed from ERA5 and SARAH-3 data. Only a share of the renewable potentials may be built (Table~\ref{tab:sys_params}). The \texttt{PyPSA-Eur} configuration file is provided in the supplementary material, so that these data can be reproduced.

The techno-economic parameters of generation and storage (Table~\ref{tab:tech_params}) refer to 2030 and are taken from the file \texttt{costs\_2030.csv} of the technology-data project, release v0.15.0 \citep{technologu_data_2024}. Monetary values are in \euro\textsubscript{2025}, and the primary source of each value, mostly the Danish Energy Agency (DEA) \citep{dea_technology_data}, is listed in the file.

Investment costs are annualized over the lifetime at the discount rate in Table~\ref{tab:sys_params} and added to the fixed O\&M costs, and the transmission costs from \texttt{PyPSA-Eur} are computed in the same way. Variable costs of dispatchable generators combine variable O\&M with fuel costs divided by electrical efficiency, and their emission factors are the CO\textsubscript{2} intensity of the fuel divided by electrical efficiency. Biogas power is modeled as a biogas plant combined with a CCGT unit and is treated as carbon-neutral.

\begin{table}[htbp]
    \small
    \centering
    \caption{General system parameters.}
    \label{tab:sys_params}
    \begin{tabular}{@{} l S[table-format=3.3, group-separator={,}] l @{}}
        \toprule
        \textbf{Parameter} & {\textbf{Value}} & \textbf{Unit} \\
        \midrule
        Discount rate                                   & 7      & \% \\
        Clustering weight of Germany                    & 0.82  & {--} \\
        Clustering weight of each neighboring country   & 0.018 & {--} \\
        Load shedding cost & 1 & \si{\Currency\per\kWh} \\
        Involuntary loss of load cost & 100 & \si{\Currency\per\kWh} \\
        System carbon intensity cap & 0.005 & \si{\TonCOtwo\per\MWh}  \\
        Max demand curtailment ratio & 5.0 & \si{\percent} \\
        Time-step duration & 25 & \si{\hour} \\
        Usable share of onshore wind potential & 20 & \si{\percent} \\
        Usable share of offshore wind potential & 50 & \si{\percent} \\
        Usable share of solar PV potential & 15 & \si{\percent} \\
        \bottomrule
    \end{tabular}
\end{table}

\begin{table}[htbp]
    \small
    \centering
    \begin{threeparttable}
        \caption{Techno-economic parameters of generation and storage technologies for 2030 \citep{technologu_data_2024}. Monetary values are in \euro\textsubscript{2025}, and all values of generators refer to their electrical output, so their efficiencies are already included.}
        \label{tab:tech_params}

        \setlength{\tabcolsep}{5pt}   
        \begin{tabular}{@{} l S[table-format=3.0] S[table-format=5.2] S[table-format=3.2] S[table-format=3.2] S[table-format=1.3] S[table-format=1.1] S[table-format=3.2] @{}}
            \toprule
            \textbf{Technology} & {\textbf{Lifetime}} & {\textbf{CAPEX}} & {\textbf{Fixed O\&M}} & {\textbf{Var.\ cost\tnote{a}}} & {\textbf{Emissions\tnote{b}}} & {\textbf{Capacity}} & {\textbf{Annualized}\tnote{d}} \\
            & {(\si{\year})} & {(\si{\Currency\per\kW})} & {(\si{\Currency\per\kW\per\year})} & {(\si{\Currency\per\MWh})} & {(\si{\TonCOtwo\per\MWh})} & {\textbf{credit}\tnote{c}} & {(\si{\Currency\per\kW\per\year})} \\
            \midrule
            Onshore wind                            & 30  & 1383.31  & 16.83  & 1.80   & 0.000 & 0.1  & 128.31 \\
            Offshore wind                           &     &          &        & 0.03   & 0.000 & 0.1  & 295.98 \\
            \quad Turbine                           & 30  & 2114.99  & 49.04  &        &       &      &        \\
            \quad AC substation                     & 40  & 850.00   & 12.75  &        &       &      &        \\
            Solar PV (utility-scale)                & 40  & 482.48   & 11.94  & 0.01   & 0.000 & 0.0  & 48.14  \\
            CCGT (gas)                              & 25  & 1108.72  & 37.14  & 54.60  & 0.341 & 1.0  & 132.27 \\
            Biogas                                  &     &          &        & 141.49 & 0.000 & 1.0  & 488.77 \\
            \quad Biogas plant\tnote{e}             & 20  & 2070.67  & 161.03 &        &       &      &        \\
            \quad CCGT power unit                   & 25  & 1108.72  & 37.14  &        &       &      &        \\
            Nuclear                                 & 40  & 10805.70 & 137.23 & 27.32  & 0.000 & 1.0  & 947.76 \\
            \midrule
            Electrolyzer\tnote{f}                   & 25  & 1886.00  & 75.44  & {--}   & 0.000 & {--} & 237.28 \\
            Fuel cell                               & 10  & 1469.38  & 73.47  & {--}   & 0.000 & {--} & 282.68 \\
            H\textsubscript{2} storage \tnote{g}    & 100 & 2.67     & 0.00   & {--}   & 0.000 & {--} & 0.19   \\
            \bottomrule
        \end{tabular}

        \begin{tablenotes}[flushleft]
            \footnotesize
            \item[a] Variable O\&M plus thermal fuel cost divided by electrical efficiency. Fuel prices (\si{\Currency\per\MWh}) and efficiencies are: natural gas (28.42, 58\%), biogas (78.81, 58\% via CCGT), and nuclear (7.45, 32.6\%).
            \item[b] CO\textsubscript{2} intensity of natural gas (\SI{0.198}{\TonCOtwo\per\MWh} thermal) divided by the CCGT efficiency.
            \item[c] Used in the reserve-margin constraint, following \citet{47_pypsa_user_meeting_2024}.
            \item[d] Annuity of the CAPEX over the lifetime at 7\% discount rate, plus fixed O\&M, summed over the components.
            \item[e] CAPEX and fixed O\&M are derived by dividing the technology data by the CCGT efficiency.
            \item[f] Electrolyzer CAPEX is per \si{\kW} electricity input (62.17\% efficiency); fuel-cell CAPEX is per \si{\kW} electricity output (50\% efficiency).
            \item[g] CAPEX per \si{\kWh} of hydrogen, and fixed O\&M and annualized cost per \si{\kWh} per year.
        \end{tablenotes}

    \end{threeparttable}
\end{table}
\section{Detailed Results} \label{app:detailed_results}

\begin{table}[htbp]
    \centering
    \small
    \setlength{\tabcolsep}{4pt}
    \begin{threeparttable}
    \caption{Investment decisions of the deterministic, reserve-margin, and robust capacity plans.\tnote{a}}
    \label{tab:results_investment}
    \begin{tabular}{@{} l l *{5}{S[table-format=5.0, group-separator={,}, group-minimum-digits=4]} @{}}
        \toprule
        & & & & \multicolumn{3}{c}{\textbf{Robust}} \\
        \cmidrule(lr){5-7}
        \textbf{Metric} & \textbf{Unit} & {\textbf{Deterministic}} & {\textbf{Reserve Margin}} & {\textbf{Case A}} & {\textbf{Case B}} & {\textbf{Case C}} \\
        \midrule
        Total Generation Capacity & \si{\GW} & 997 & 1068 & 984 & 991 & 968 \\
        \quad On-wind Capacity & \si{\GW} & 471 & 458 & 459 & 465 & 440 \\
        \quad Off-wind Capacity & \si{\GW} & 70 & 69 & 70 & 70 & 70 \\
        \quad Solar Capacity & \si{\GW} & 291 & 279 & 290 & 293 & 293 \\
        \quad Gas Capacity & \si{\GW} & 65 & 144 & 63 & 63 & 61 \\
        \quad Biogas Capacity & \si{\GW} & 63 & 82 & 65 & 63 & 63 \\
        \quad Nuclear Capacity & \si{\GW} & 36 & 36 & 38 & 39 & 41 \\
        \midrule
        Electrolyzer & \si{\GW} & 9 & 0 & 9 & 11 & 14 \\
        Fuel Cell & \si{\GW} & 17 & 0 & 16 & 17 & 19 \\
        H\textsubscript{2} Storage & \si{\GWh} & 9536 & 0 & 9556 & 9322 & 10122 \\
        \midrule
        Total Transmission Capacity & \si{\GW} & 494 & 502 & 493 & 503 & 522 \\
        \quad \% of Max Potential & \si{\percent} & 33 & 33 & 33 & 33 & 35 \\
        \quad Increase over Existing Grid\tnote{b} & \si{\percent} & 55 & 57 & 54 & 57 & 64 \\
        \bottomrule
    \end{tabular}
    \begin{tablenotes}[flushleft]
        \footnotesize
        \item[a] The reserve-margin model uses $\RESERVEMARGIN = 1.1$. Robust plans are the incumbent solutions from Algorithm~\ref{alg:Benders_Decomposition} for outage Cases~A, B, and C. Figure~\ref{fig:all_plans} shows the spatial distribution of these capacity mixes.
        \item[b] Relative to the existing transmission capacity of \SI{319}{\GW}.
    \end{tablenotes}
    \end{threeparttable}
\end{table}

\begin{landscape}
    \begin{table}[htbp]
        \centering
        \caption{Bounds, computation times, and simulated cost and loss-of-load (LoL) statistics of the deterministic, reserve-margin, and robust capacity plans across test cases\textsuperscript{a}. All cost and LoL statistics are computed over $\noSimScenarios = 1000$ out-of-sample simulations per case.}
        \label{tab:results_cases}
        \resizebox{\linewidth}{!}{%
        \begin{threeparttable}
        \begin{tabular}{@{} l l *{9}{S[table-format=4.2, group-separator={,}, group-minimum-digits=4]} @{}}
            \toprule
            & & \multicolumn{3}{c}{\textbf{Case A}} & \multicolumn{3}{c}{\textbf{Case B}} & \multicolumn{3}{c}{\textbf{Case C}} \\
            & & \multicolumn{3}{c}{MTTR = \SI{12}{\hour}, $\FailRateBase$ = 0.5 failures/100~km/year} & \multicolumn{3}{c}{MTTR = \SI{24}{\hour}, $\FailRateBase$ = 1.0 failures/100~km/year} & \multicolumn{3}{c}{MTTR = \SI{48}{\hour}, $\FailRateBase$ = 2.0 failures/100~km/year} \\
            \cmidrule(lr){3-5} \cmidrule(lr){6-8} \cmidrule(lr){9-11}
            \textbf{Metric} & \textbf{Unit} & {\textbf{Deterministic}} & {\textbf{Reserve Margin}} & {\textbf{Robust}}
            & {\textbf{Deterministic}} & {\textbf{Reserve Margin}} & {\textbf{Robust}}
            & {\textbf{Deterministic}} & {\textbf{Reserve Margin}} & {\textbf{Robust}} \\
            \midrule
            
            System Cost Lower Bound ($\LB$) & \si{B\Currency} &  &  & 213.6 &  &  & 213.9 &  &  & 215.0 \\
            System Cost Upper Bound ($\UB$) & \si{B\Currency} &  &  & 214.2 &  &  & 215.2 &  &  & 217.1 \\
            Relative Gap & \si{\percent} &  &  & 0.28 &  &  & 0.61 &  &  & 0.94 \\
            Computation Time\tnote{b} & \si{\hour} & 0.15 & 0.13 & 17.85 & 0.15 & 0.12 & 17.88 & 0.16 & 0.11 & 14.73 \\
            \midrule \addlinespace
    
            Highest Adjusted System Cost\tnote{c} & \si{B\Currency} & 213.9 & 224.5 & 213.8 & 215.4 & 225.9 & 215.0 & 218.4 & 228.1 & 216.0 \\
            Mean Adjusted System Cost & \si{B\Currency} & 213.1 & 223.7 & 213.3 & 213.5 & 224.0 & 213.9 & 215.2 & 225.3 & 215.4 \\
            Lowest Adjusted System Cost & \si{B\Currency} & 213.0 & 223.6 & 213.3 & 213.0 & 223.6 & 213.7 & 213.3 & 223.9 & 214.9 \\
            Std. Adjusted System Cost & \si{B\Currency} & 0.2 & 0.1 & 0.1 & 0.4 & 0.3 & 0.1 & 0.8 & 0.6 & 0.2 \\
            \midrule \addlinespace
    
            Annualized Investment and O\&M Cost & \si{B\Currency} & 183.9 & 192.4 & 184.0 & 183.9 & 192.4 & 185.6 & 183.9 & 192.4 & 186.3 \\
            Mean Operational Cost & \si{B\Currency} & 28.6 & 31.2 & 28.7 & 28.7 & 31.4 & 27.7 & 29.5 & 32.2 & 28.5 \\
            Mean Load Shedding Cost & \si{B\Currency} & 0.6 & 0.0 & 0.6 & 0.8 & 0.2 & 0.6 & 1.7 & 0.7 & 0.6 \\
            \midrule \addlinespace
    
            Total Load & \si{\TWh} & 1701 & 1701 & 1701 & 1701 & 1701 & 1701 & 1701 & 1701 & 1701 \\
            Highest Loss of Load & \si{\GWh} & 408 & 0 & 0 & 701 & 0 & 0 & 1371 & 6 & 2 \\
            Mean Loss of Load & \si{\GWh} & 3 & 0 & 0 & 18 & 0 & 0 & 79 & 0 & 0 \\
            Lowest Loss of Load & \si{\GWh} & 0 & 0 & 0 & 0 & 0 & 0 & 0 & 0 & 0 \\
            Std. Loss of Load & \si{\GWh} & 27 & 0 & 0 & 77 & 0 & 0 & 172 & 0 & 0 \\
            \bottomrule
        \end{tabular}%
        
        \begin{tablenotes}
            \item[a] For the reserve-margin model, the reported solution corresponds to the reserve margin level $\Gamma = 1.1$. The reported robust plans (the incumbents of Algorithm~\ref{alg:Benders_Decomposition}) were found at iterations 108, 139, and 242 for Cases~A, B, and C, respectively. $\LB$ is the master-problem lower bound at the end of Phase~I (300 iterations). $\UB$ is the simulation-based Cantelli bound on the VaR of the total system cost at probability level $\alpha = 0.99$, computed in Phase~II from the 1000 out-of-sample simulations, and the relative gap is $\Gap = (\UB - \LB) / \UB \times 100$.
            \item[b] Wall-clock time of the optimization plus the out-of-sample simulation.
            \item[c] The adjusted system cost re-prices involuntary loss of load (LoL) from its punitive simulation rate of \SI{100}{\Currency\per\kWh} to the same \SI{1}{\Currency\per\kWh} penalty applied to voluntary load shedding, as defined in Eq.~\eqref{eq:adjusted_cost}. This adjustment eliminates the distortion caused by the high LoL penalty, while LoL itself is reported separately in \si{\GWh} in this table.
        \end{tablenotes}
        \end{threeparttable}}
    \end{table}
\end{landscape}

\begin{figure}[p]
    \centering
    \includegraphics[width=\textwidth]{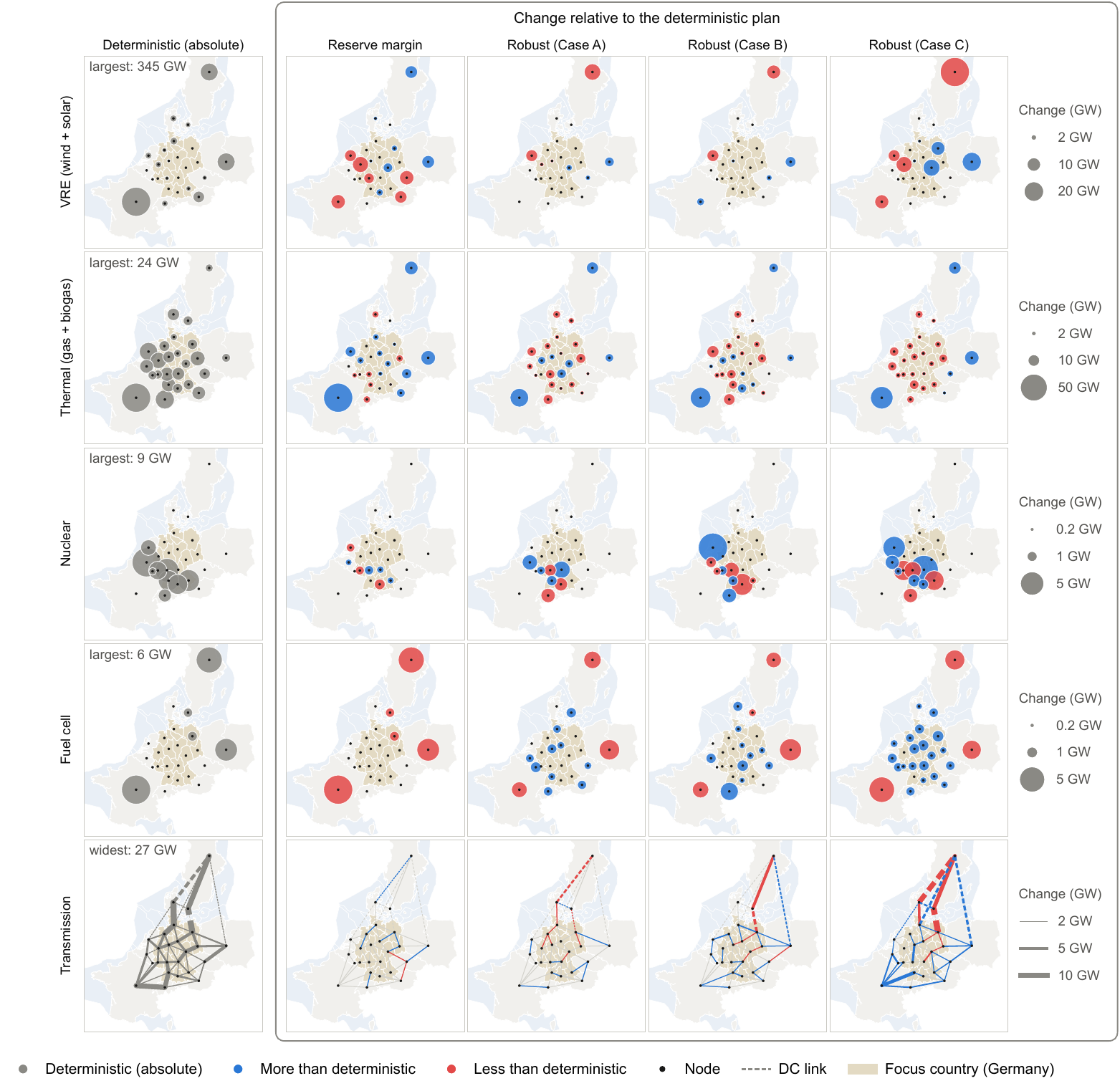}
    \caption{Spatial distribution of capacity in the deterministic plan (first column) and changes relative to it in the reserve-margin plan ($\Gamma = 1.1$) and the robust plans for Cases~A, B, and~C (framed columns). Rows show variable renewable energy (VRE) capacity (onshore wind, offshore wind, and solar), thermal capacity (gas and biogas), nuclear capacity, fuel-cell capacity, and transmission capacity. In the first column, circle area is proportional to the installed capacity at each node, and the width of each drawn line to the capacity of that line. In the framed columns, circle area and line width are proportional to the change in capacity, with blue denoting more and red less capacity than in the deterministic plan; within each row, all framed panels share the same scale. Changes smaller than 0.3~\si{\GW} at nodes and 0.5~\si{\GW} on lines are not shown.}
    \label{fig:all_plans}
\end{figure}

\end{document}